\documentclass[mn,a4paper,fleqn]{w-art}
\usepackage{amsfonts}
\usepackage{times,cite,w-thm}

\theoremstyle{plain}

\theoremstyle{definition}

\usepackage[]{graphicx}

\numberwithin{equation}{section}

\newcommand{\eps}{\varepsilon}

\newcommand{\norm}[1]{\big\Vert#1\big\Vert}

\newcommand{\abs}[1]{\left\vert#1\right\vert}
\newcommand{\set}[1]{\left\{\,#1\,\right\}}

\newcommand{\inner}[1]{\left(#1\right)}
\newcommand{\biginner}[1]{\big(#1\big)}
\newcommand{\comi}[1]{\left<#1\right>}
\newcommand{\com}[1]{\left[#1\right]}
\newcommand{\reff}[1]{(\ref{#1})}

\begin{document}


\title{Global Hypoellipticity and Compactness of Resolvent for Fokker-Planck Operator}

%\author[F. Author]{First Author\inst{1,}%
%  \footnote{Corresponding author\quad E-mail:~\textsf{x.y@xxx.yyy.zz},
%            Phone: +00\,999\,999\,999,
%            Fax: +00\,999\,999\,999}}
%\address[\inst{1}]{First address}

\author{Wei-Xi Li}
\address{School of Mathematics and Statistics,
Wuhan University, 430072 Wuhan, China}
\address{Institut de Math\'{e}matiques de Jussieu, Universit\'{e} Paris 6,
75013 Paris, France}
\address{E-mail: wei-xi.li@whu.edu.cn}

\begin{abstract}
  {\bf Abstract}~In this paper we study the Fokker-Planck operator  with potential
  $V(x),$ and
  analyze some kind of conditions imposed on
  the potential to ensure the validity of global hypoelliptic
  estimates (see Theorem \ref{+Hypo}).  As a consequence, we obtain the compactness of resolvent of
  the Fokker-Planck operator if either the Witten Laplacian on 0-forms has a compact resolvent or some additional assumption on the
  behavior of the potential at infinity is fulfilled. This work improves the previous results
  of H\'{e}rau-Nier \cite{HerauNier04} and Helffer-Nier \cite{HelfferNier05}, by obtaining a better
  global hypoelliptic estimate under weaker assumptions
  on the potential.
\end{abstract}
\pagestyle{plain} \maketitle

\section{Introduction and main results}

In this work we consider the Fokker-Planck operator
\begin{eqnarray}\label{FP++}
  P=y\cdot\partial_x-\partial_x V(x)\cdot
  \partial_y-\triangle_y+\frac{\abs y^2}{4}-\frac{n}{2},\quad
  (x,y)\in\mathbb{R}^{2n}
\end{eqnarray}
where $x$ denotes the space variable and $y$ denotes the velocity
variable, and $V(x)$ is a potential defined in the whole space
$\mathbb{R}_x^n$. There have been extensive works concerned with the
operator $P$, with various techniques  from different fields such as
partial differential equation, spectral theory and statistical
physics. In this paper we will focus on analyzing some  kind of
conditions imposed on the potential $V(x),$ so that the
Fokker-Planck operator $P$ admits a global hypoelliptic estimate and
has a compact resolvent. This problem is linked closely with the
trend to equilibrium for the Fokker-Planck operator, and has been
studied by Desvillettes-Villani, Helffer-Nier, H\'{e}rau-Nier and
some other authors (see
\cite{DesvVillani01,HelfferNier05,HerauNier04} and the references
therein). It is believed that the global estimate and the
compactness of resolvent are related to the properties of the
potential $V(x).$ In the particular case of quadratic potential, the
theory is well developed. As far as general potential is concerned,
different kind of assumptions on $V(x)$ had been explored firstly by
H\'{e}rau-Nier \cite{HerauNier04} and then generalized by
Helffer-Nier \cite{HelfferNier05}. This work is motivated  by the
previous works of H\'{e}rau-Nier and Helffer-Nier, and can be seen
as an improvement of their  results. Our main result is the
following.

\begin{theorem}\label{+Hypo}
Let $V(x)\in C^2(\mathbb{R}^n)$ be a real-valued function satisfying
that
\begin{eqnarray}\label{+HyP}
   \forall~\abs\alpha=2,~\exists~C_\alpha>0,\quad \abs{\partial_x ^\alpha V(x)}\leq C_\alpha \inner{1+\abs{\partial_x V(x)}^2}^{s\over2}
    ~~\,\,{\rm with~~s<\frac{4}{3}} .
\end{eqnarray}
Then there is a constant $C,$ such that for any $u\in
 C_0^\infty\big(\mathbb{R}^{2n}\big)$ one has
  \begin{eqnarray}\label{+A1}
   \norm{ \abs{\partial_x  V(x)}^{2\over3}u}_{L^2}
   \leq C\set{\,
   \norm{P u}_{L^2}+\norm{u}_{L^2}\,},
\end{eqnarray}
and
  \begin{eqnarray}\label{+++A1}
   \norm{\inner{1-\triangle_x}^{\frac{\delta}{2}}u}_{L^2}
   +\norm{ \inner{1-\triangle_y+\abs y^2}^{1\over2}u}_{L^2}
   \leq C\Big\{\,
   \norm{P u}_{L^2}+\norm{u}_{L^2}\Big\},
\end{eqnarray}
where $\delta$  equals to $\frac{2}{3}$ if $s\leq \frac{2}{3}$,
$\frac{4}{3}-s$ if ${2\over3}<s\leq \frac{10}{9},$ and
$\frac{2}{3}-{s\over2}$ if ${{10}\over9}< s<\frac{4}{3}$. Here and
throughout the paper we will use $\Vert \cdot\Vert_{L^2}$ to denote
the norm of the complex Hilbert space
$L^2\big(\mathbb{R}^{2n}\big),$ and denote by $
C_0^\infty\inner{\mathbb{R}^{2n}}$ the set of smooth compactly
supported functions.
\end{theorem}

\begin{remark}
  In particular, if the assumption \reff{+HyP} is fulfilled with
  $s=\frac{2}{3}$, then we have the following hypoelliptic estimate which seems
  to be optimal:
  \[
   \forall~u\in
 C_0^\infty\big(\mathbb{R}^{2n}\big),\quad \norm{ \abs{\partial_x
V(x)}^{2\over3}u}_{L^2}+\norm{\inner{1-\triangle_x}^{\frac{1}{3}}u}_{L^2}
   \leq C\Big\{
   \norm{P u}_{L^2}+\norm{u}_{L^2}\Big\}.
  \]
  Moreover one can deduce from the above estimate a better regularity in the velocity
  variable $y,$ that is,
  \[
   \forall~u\in
 C_0^\infty\big(\mathbb{R}^{2n}\big),\quad \norm{ \inner{1-\triangle_y+\abs
y^2}u}_{L^2}\leq C\Big\{\,
   \norm{P u}_{L^2}+\norm{u}_{L^2}\Big\}.
  \]
  This can be seen in Proposition \ref{855+++} in the next section.
\end{remark}

\begin{corollary}
The Fokker-Planck operator $P$ has a compact resolvent if  the
potential $V(x)$ satisfies  \reff{+HyP} and that $ \lim\limits_{\abs
x\rightarrow+\infty}\abs{\partial_xV(x)}=+\infty.$
\end{corollary}

To analyze the compactness of resolvent of the operator $P,$ the
hypoellipticity techniques are an efficient tool, one of which is
referred to Kohn's method \cite{Kohn78} and another is based on
nilpotent Lie groups (see \cite{Hel-Nou,RothschildStein76}). Kohn's
method had been used by H\'{e}rau-Nier \cite{HerauNier04} to study
such a potential $V(x)$ that behaves at infinity as a high-degree
homogeneous function.  More precisely, if $V(x)$ satisfies that for
some $C,M\geq1,$
\begin{eqnarray}\label{7235}
  \frac{1}{C}\comi{x}^{2M-1}\leq \inner{1+\abs
  {\partial_xV(x)}^2}^{1\over2}~~ {\rm and}~~
  \forall~\abs\gamma\geq0, ~\abs{\partial_x^\gamma V(x)}\leq C_\gamma \comi{x}^{2M-\abs\gamma}
  ,
\end{eqnarray}
where $\comi{x}=\inner{1+\abs{x}^2}^{1\over2}$, then H\'{e}rau-Nier
established the following  isotropic hypoelliptic estimate, by use
of the global pseudo-differential calculus,
\begin{eqnarray}\label{7231}
  \norm{ \Lambda_{x,y}^{1\over 4} u}_{L^2}
   \leq C\set{
   \norm{P u}_{L^2}+\norm{u}_{L^2}}
\end{eqnarray}
with
$\Lambda_{x,y}=\inner{1-\triangle_x-\triangle_y+\frac{1}{2}\abs{V(x)}^2+\frac{1}{2}\abs
{y}^2}^{1\over2}.$  By developing the approach of H\'{e}rau-Nier,
Helffer-Nier \cite{HelfferNier05} obtained the same estimate as
above for more general $V(x)$ which satisfies that for some $C,
k\geq1,$
\begin{eqnarray}\label{7232}
 \frac{1}{C}\comi{x}^{\frac{1}{k}}\leq \Big(1+|\partial_x  V(x)|^2\Big)^{1\over2}\leq
  C\comi{x}^{k}~\textrm{and}~\forall\abs\gamma\geq 0, ~\abs{\partial_x ^\gamma V(x)}\leq  C_\gamma
  \big(1+\abs{\partial_x  V(x)}^2\big)^{1\over2}.
\end{eqnarray}
As for the Kohn's proof for the hypoellipticity, the exponent
${1\over4}$ in \reff{7231} is not optimal.  A better exponent, which
seems to be ${2\over3}$  as seen in  \cite{RothschildStein76}, can
be obtained  via explicit method in the particular case  when $V(x)$
is a non-degenerate quadratic form.   Moreover Helffer-Nier
\cite{HelfferNier05} studied such a $V(x)$ that satisfies
\begin{eqnarray}\label{7261}
  \forall~\abs\alpha=2,\quad \abs{\partial_x^\alpha V(x)}\leq C_\alpha
  \inner{1+\abs{\partial_x  V(x)}^2}^{\frac{1-\rho}{2}}\quad {\rm
  with}~~\,
  \rho>\frac{1}{3},
\end{eqnarray}
and obtained the estimate
\begin{eqnarray}\label{A+B}
   \norm{ \abs{\partial_x  V(x)}^{2\over3}u}_{L^2}
   \leq C\set{\,
   \norm{P u}_{L^2}+\norm{u}_{L^2}\,}.
\end{eqnarray}
This generalized the quadratic potential case, and their main tool
is the nilpotent technique that initiated by
\cite{RothschildStein76} and then developed by \cite{Hel-Nou}.
Although the estimate \reff{A+B} is better, the condition
\reff{7261} is stronger than \reff{7232} for the second derivatives,
and comparing with \reff{7231}, we see that in \reff{A+B} some
information on the Sobolev regularity in $x$ is missing. In
\reff{+HyP} we get rid of the assumptions on the behavior of
$\partial_xV(x)$ at infinity. This generalizes the conditions
\reff{7235} and \reff{7232}. Moreover, the exponent in \reff{+A1} is
${2\over3}$, better than ${1\over4}$ established in \reff{7231}.
Besides, we have relaxed the condition \reff{7261} by allowing the
number $\rho$ there to take values in the interval
$]-{1\over3},~+\infty[.$  As seen in the proof presented in Section
\ref{sec3}, our approach is direct, which seems simpler for it
doesn't touch neither complicated nilpotent group techniques nor
pseudo-differential calculus.

%Now we analyze the compactness of resolvent of the Fokker-Planck
%operator $P.$  As an immediate consequence of \reff{+A1} and
%\reff{+++A1}, we see the form domain of the  operator $P$ is
%included in the space $\set{u\in H^\delta\big(\mathbb{R}^{2n}\big);~
%\inner{\abs {\partial_xV(x)}+\abs y}^{\delta} u\in
%L^2\big(\mathbb{R}^{2n}\big)}$ with $\delta>0,$ which is a compactly
%embedded in $L^2\big(\mathbb{R}^{2n}\big)$  if
%\begin{eqnarray}\lab{7238}
%  \lim\limits_{\abs
%  x\rightarrow+\infty}\abs{\partial_xV(x)}=+\infty.
%\end{eqnarray}
%We state this result as the following corollary.

% We remark that \reff{7238} is not a necessary condition for the
%  compactness of resolvent of the operator $P.$ This can be seen in
%  the case
%  when the dimension $n=2$ and $\abs{\partial_xV(x)}^2\approx x_1^2x_2^2.$
%  In such a case we can verify that the following inequalities hold
%  in the sense of operators:
%  \[
%    -\triangle_x+\abs{\partial_xV(x)}^2\geq
%    -\triangle_x+x_1^2x_2^2\geq
%    \frac{1}{2}\inner{\abs{x_1}+\abs{x_2}}.
%  \]
% As a result, if \reff{+A1} and \reff{+++A1} hold  then the form domain of the Fokker-Planck operator $P$ is included
%in the space $\set{u(x_1,x_2)\in L^2\big(\mathbb{R}^{2}\big);~
%\inner{\abs {x_1}+\abs{x_2}}^{\delta\over2} u\in
%L^2\big(\mathbb{R}^{2}\big)},$ which is compact in
%$L^2\inner{\mathbb{R}^{2}}$, and hence $P$ has a compact resolvent.

 Another direction to
get the compact resolvent is to analyze the relationship between $P$
and the Witten Laplace operator $\triangle_{V/2}^{(0)}$ defined by
\[
  \triangle_{V/2}^{(0)}=
  -\triangle_x+\frac{1}{4}\abs{\partial_xV(x)}^2-\frac{1}{2}\triangle_x
  V(x).
\]
In \cite{HelfferNier05}, Helffer-Nier stated a conjecture which says
that the Fokker-Planck operator $P$ has a compact resolvent if and
only if the Witten Laplacian $\triangle_{V/2}^{(0)}$ has a compact
resolvent. The necessity part is well-known, and under rather weak
assumptions on the potential $V,$ saying $V\in
C^\infty\big(\mathbb{R}^{2n}\big)$ for instance, if the
Fokker-Planck operator $P$ has a compact resolvent then the Witten
Laplacian $\triangle_{V/2}^{(0)}$ has a compact resolvent.  The
reverse implication still remains open, and some partial answers
have been obtained by \cite{HelfferNier05,HerauNier04}. For example,
suppose $V\in C^\infty\big(\mathbb{R}^{2n}\big)$ such that
\begin{eqnarray*}
  \forall\abs\gamma\geq 0,~~\forall~x\in\mathbb{R}^{2n},\quad \abs{\partial_x ^\gamma V(x)}\leq  C_\gamma
  \inner{1+\abs{\partial_x  V(x)}^2}^{1\over2},
\end{eqnarray*}
\begin{eqnarray*}
  \exists~M,~C>1,~~\forall~x\in\mathbb{R}^{2n},\quad \abs{\partial_x  V(x)}\leq
  C\comi{x}^{M},
\end{eqnarray*}
and
\begin{eqnarray*}
  \exists~\kappa>0,~~\forall~ \abs\alpha=2,~~\forall~x\in\mathbb{R}^{2n},
  \quad\abs{\partial_x ^\alpha V(x)}\leq
  C_\alpha
  \inner{1+\abs{\partial_x  V(x)}^2}^{1\over2}\comi{x}^{-\kappa}.
\end{eqnarray*}
Then the operator $P$ has a compact resolvent if the Witten Laplace
operator $\triangle_{V/2}^{(0)}$ has a compact resolvent (see
Corollary 5.10 of \cite{HelfferNier05}). Due to Theorem \ref{+Hypo},
we can generalize the previous results as follows.

\begin{corollary}\label{thm+2}
  Let $V(x)$ satisfy the condition \reff{+HyP}.
  Then the Fokker-Planck operator $P$ has a compact resolvent if the Witten
  Laplacian $\triangle_{V/2}^{(0)}$ has a compact resolvent.
\end{corollary}

The paper is organized as follow. In the next section we introduce
some notations used throughout the paper, and then present some
regularity results on the velocity variable $y$. Since the proof of
Theorem \ref{+Hypo} is quite lengthy, we divide it into two parts
and proceed to handle them in Section \ref{sec3} and Section
\ref{sec4}. The proof of  Corollary \ref{thm+2} will be presented in
Section \ref{sec5}.

\section{Notation and regularity in velocity variable}\label{sec2}
We  firstly list some notations used throughout the paper in
Subsection \ref{sec21}, and then establish the regularity in the
velocity variable $y$ in Subsection \ref{sec22}.  This will give the
desired estimate on the second term on the left of \reff{+++A1}.

\subsection{Notation}\label{sec21}
Throughout the paper we denote by $(\xi,\eta)$ the dual variables of
$(x,y),$  and denote by $\comi{\cdot,~\cdot}_{L^2}$ the inner
product of the complex Hilbert space $L^2\big(\mathbb{R}^{2n}\big).$
Set
\[
D_{x_j}=-i\,\partial_{x_j}, \quad D_{y_j}=-i\,\partial_{y_j}~~ {\rm
and}~~D_x=(D_{x_1},\cdots, D_{x_n}),\quad D_y=(D_{y_1},\cdots,
D_{y_n}).
\]
Let $\Lambda_y$ be the  operator given by
\[\Lambda_{y}=\inner{1+{1\over2}\abs{y}^2-\triangle_y}^{1\over2}.\]
Observing $\abs{\partial_xV(x)}$ is only continuous, we have to
replace it sometimes by the equivalent $C^1$ function $f(x)$ given
by
\[
f(x)=\inner{1+\abs{\partial_x V(x)}^2}^{1\over2}.
\]
Denoting $Q=y \cdot D_x-\partial_x V(x)\cdot
 D_y$ and $L_j=\partial_{y_j}+\frac{y_j}{2},j=1,\cdots n,$ we can
 write the operator $P$ given in \reff{FP++} as
\begin{equation}\label{+FP+}
  P=iQ+\sum_{j=1}^nL_j^*L_j.
\end{equation}

\subsection{Regularity  in the velocity variable}\label{sec22} In view of the expression \reff{+FP+}, we see
that the required estimate on the  term $\norm{\Lambda_yu}_{L^2}$ is
easy to get, without any assumption on the potential $V(x)$. Indeed,
As a result of \reff{+FP+}, we have
\begin{eqnarray}\label{yva}
  \forall~u\in  C_0^\infty\big(\mathbb{R}^{2n}\big),\quad\sum_{j=1}^n\norm{L_j
  u}_{L^2}^2\leq {\rm Re}\comi{Pu,~u}_{L^2},
\end{eqnarray}
from which one can deduce that
\begin{eqnarray}\label{velo}
   \forall u\in  C_0^\infty\big(\mathbb{R}^{2n}\big),\quad\norm{
   \Lambda_yu}_{L^2}^2
   \leq C\Big\{
   \abs{\comi{P u,~u}_{L^2}}+\norm{u}_{L^2}^2\Big\}.
\end{eqnarray}
This gives the desired estimate on the second term on the left of
\reff{+++A1}.

For constant potential, i.e., $\partial_xV(x)=0,$ starting from the
regularity  in $x$, we can derive a better Sobolev exponent, which
is known to be $2,$ for the regularity in $y$ variable (see for
instant \cite{Bouchut02}). When general potential is involved, we
have the following estimate.

\begin{proposition}\label{855+++}
  There exists a constant $C$ such that  for any  $u\in  C_0^\infty\big(\mathbb{R}^{2n}\big),$
  \begin{eqnarray}\label{twore}
    \norm{\Lambda_y^2u}_{L^2}
    \leq C\Big\{\norm{\abs{\partial_xV(x)}^{2\over3} u}_{L^2}+\norm{\inner{1-\triangle_x}^{1\over3}u}_{L^2}
    + \norm{P u}_{L^2} \Big\},
  \end{eqnarray}
  or equivalently,
  \begin{eqnarray}\label{s2+}
   \sum_{j=1}^n\norm{L_jL_j^*u}_{L^2} \leq C\set{\norm{\abs{\partial_{x}V(x)}^{2\over3}~u}_{L^2}
   +\norm{\inner{1-\triangle_x}^{1\over3}~u}_{L^2}
   +\norm{Pu}_{L^2} }.
\end{eqnarray}
\end{proposition}

\noindent\emph{Proof.} In this proof we show \reff{s2+}. Using
\reff{yva} gives
\begin{eqnarray*}
  \norm{L_jL_j^*u}_{L^2}^2
  &\leq& {\rm Re}\comi{PL_j^*u,~L_j^*
  u}_{L^2}\\
  &=&{\rm Re}\comi{[P,~L_j^*]u,~L_j^*
  u}_{L^2}+ {\rm Re}\comi{Pu,~L_jL_j^*
  u}_{L^2}\\
  &\leq&{\rm Re}\comi{[P,~L_j^*]u,~L_j^*
  u}_{L^2}+{1\over2}
  \norm{L_jL_j^*u}_{L^2}^2+2\norm{Pu}_{L^2}^2.
\end{eqnarray*}
Hence
\[
  \norm{L_jL_j^*u}_{L^2}^2
  \leq 2\abs{\comi{[P,~L_j^*]u,~L_j^*
  u}_{L^2}}+4\norm{Pu}_{L^2}^2.
\]
Now assume the following estimate holds, for any $\eps>0$,
\begin{eqnarray}\label{0909281}
   \abs{\comi{[P, L_j^*]u, L_j^*u}_{L^2}}
   \leq \eps \norm{L_jL_j^*u}_{L^2}^2 +
   C_\eps\Big\{\norm{|\partial_{x}V|^{2\over3}u}_{L^2}^2
   +\norm{(1-\triangle_x)^{1\over3}u}_{L^2}^2
   +\norm{Pu}_{L^2}^2\Big\}.
\end{eqnarray}
Then combining the above two inequalities and then letting $\eps$
small enough, we get the desired estimate \reff{s2+}.  In order to
show \reff{0909281}, we make use of the following commutation
relations satisfied by $iQ, L_j, L_k^*, j,k=1,2,\cdots,n,$
\[
  [iQ,~L_j^*]=-\frac{1}{2}\partial_{x_j}
   V(x)+\partial_{x_j},\quad [L_j,~L_k]=[L_j^*,~L_k^*]=0,\quad [L_j,~L_k^*]=\delta_{jk};
\]
this gives
\[
  [P,~L_j^*]=-\frac{1}{2}\partial_{x_j}V(x)+\partial_{x_j}+L_j^*.
\]
Then
\begin{eqnarray*}
   \abs{\comi{[P,~L_j^*]u,~L_j^*u}_{L^2}}&\leq&
   \comi{L_j^*u,~L_j^*u}_{L^2}
   +\abs{\Big<\Big(-\frac{1}{2}\partial_{x_j}
   V(x)+\partial_{x_j}\Big)u,~L_j^*u\Big>_{L^2}}\\
   &\leq& C\set{\norm{\abs{\partial_{x}V(x)}^{2\over3}~u}_{L^2}^2
   +\norm{\inner{1-\triangle_x}^{1\over3}~u}_{L^2}^2+\norm{L_ju}_{L^2}^2}\\
   &&+C\set{\norm{L_j\abs{\partial_{x}V(x)}^{1\over3}u}_{L^2}^2+
   \norm{L_j\inner{1-\triangle_x}^{1\over6}u}_{L^2}^2}.
\end{eqnarray*}
Moreover, note that
\begin{eqnarray*}
  \norm{L_j\abs{\partial_{x}V(x)}^{1\over3}u}_{L^2}^2
  &=&\comi{L_j^*L_ju,~\abs{\partial_{x}V(x)}^{2\over3}u}_{L^2}\\
  &=&\comi{L_jL_j^*u,~\abs{\partial_{x}V(x)}^{2\over3}u}_{L^2}
  -\comi{u,~\abs{\partial_{x}V(x)}^{2\over3}u}_{L^2},
\end{eqnarray*}
and hence
\begin{eqnarray*}
  \forall~\eps>0,\quad \norm{L_j\abs{\partial_{x}V(x)}^{1\over3}u}_{L^2}^2
  \leq
  \eps\norm{L_jL_j^*u}_{L^2}^2+C_\eps\set{\norm{\abs{\partial_{x}V(x)}^{2\over3}u}_{L^2}^2
  +\norm{u}_{L^2}^2}.
\end{eqnarray*}
Similarly,
\begin{eqnarray*}
  \forall~\eps>0,\quad \norm{L_j\inner{1-\triangle_x}^{1\over6}u}_{L^2}^2
  \leq
  \eps\norm{L_jL_j^*u}_{L^2}^2+C_\eps\set{\norm{\inner{1-\triangle_x}^{1\over3}u}_{L^2}^2
  +\norm{u}_{L^2}^2}.
\end{eqnarray*}
These inequalities yield \reff{0909281}. The proof of Proposition
\ref{855+++} is thus completed. \qed

\section{Proof of Theorem
\ref{+Hypo}: the first part}\label{sec3}

In this section we only show \reff{+A1} and postpone \reff{+++A1} to
the next section. Let $V(x)$ satisfy the assumption \reff{+HyP}.
Then using the notation
\[f(x)=\inner{1+\abs{\partial_xV(x)}^2}^{1\over2},\]
we have
\begin{eqnarray}\label{cf+}
  \forall~x\in\mathbb{R}^n,\quad\abs{\partial_x f}\leq C f(x)^{s}\quad{\rm with}~~s<{4\over3}.
\end{eqnarray}

The following is the main result of this section.
\begin{proposition}\label{prop+1}
  Suppose $f$ satisfies the condition \reff{cf+}. Then
  \begin{eqnarray}\label{req}
    \exists ~C>0,~~\forall~u\in C_0^\infty\big(\mathbb{R}^{2n}\big),\quad
    \norm{f(x)^{2\over3}u}_{L^2}\leq C
    \set{\norm{Pu}_{L^2}+\norm{u}_{L^2}}.
  \end{eqnarray}
\end{proposition}

\noindent\emph{Proof.} To simplify the notation, we will use the
capital letter $C$ to denote different suitable constants. Let $R\in
C^1\big(\mathbb{R}^{2n}\big)$ be a real-valued function given by
\[
   R=R(x,y)= 2f(x)^{-{2\over3}}\partial_x V(x)\cdot y.
\]
We can verify that
\begin{eqnarray*}
  \forall~u\in C_0^\infty\big(\mathbb{R}^{2n}\big),\quad \norm{R u}_{L^2}\leq
  C \norm{\abs{y}f(x)^{{1\over3}} u}_{L^2}\leq C \norm{\Lambda_{y}f(x)^{{1\over3}} u}_{L^2}.
\end{eqnarray*}
Recall $P=iQ+\sum_{j=1}^nL_j^*L_j$ with $Q=y \cdot D_x-\partial_x
V(x)\cdot
 D_y$ and $L_j=\partial_{y_j}+\frac{y_j}{2}.$  Then the above
 inequalities together with the relation
\[
  {\rm Re} \comi{Pu,~Ru}_{L^2}={\rm Re}
  \comi{iQu,~Ru}_{L^2}+
  {\rm Re} \sum_{j=1}^n\comi{L_j^*L_ju,~Ru}_{L^2}
\]
yield
\begin{eqnarray}\label{rmain}
  {\rm Re}
  \comi{iQu,~Ru}_{L^2}\leq
  \norm{Pu}_{L^2}^2+\norm{\Lambda_{y}f(x)^{{1\over3}} u}_{L^2}^2+\sum_{j=1}^n\abs{\comi{L_j^*L_ju,~Ru}_{L^2}}.
\end{eqnarray}
Next we will proceed to treat the terms on both sides of
\reff{rmain} by the following three steps.

{\bf \emph{Step I}.}~Firstly we will show that for any $\eps>0$
there exists a constant $C_\eps>0$ such that
\begin{eqnarray}\label{reqA}
  \forall~u\in C_0^\infty\big(\mathbb{R}^{2n}\big),\quad \norm{\Lambda_{y}f(x)^{{1\over3}}
  u}_{L^2}^2\leq \eps \norm{f(x)^{{2\over3}}
  u}_{L^2}^2+C_\eps\set{\norm{Pu}_{L^2}^2+\norm{u}_{L^2}^2}.
\end{eqnarray}
To confirm this, we use \reff{velo} to get
\begin{eqnarray*}
  \norm{\Lambda_{y}f(x)^{1\over3}u}_{L^2}^2
  &\leq&{\rm Re}\comi{Pf(x)^{1\over3}u,~f(x)^{1\over3}u}_{L^2}\\
  &=&{\rm Re}\comi{Pu,~f(x)^{2\over3}u}_{L^2}+
  {\rm Re}\comi{\com{P,~f(x)^{1\over3}}u,~f(x)^{1\over3}u}_{L^2}.
\end{eqnarray*}
The upper bound of the term ${\rm
Re}\comi{Pu,~f(x)^{2\over3}u}_{L^2}$ can be obtained by
Cauchy-Schwarz's inequality. Then the required estimate \reff{reqA}
will follow if the following inequality holds: for any
$\eps_1,\eps_2>0,$ there exists a constant $C_{\eps_1,\eps_2}$ such
that
\begin{eqnarray}\label{reqA1}
   \comi{\com{P,~f(x)^{1\over3}}u,~f(x)^{1\over3}u}_{L^2}\leq
  \eps_1\norm{\Lambda_{y}f(x)^{1\over3}u}_{L^2}^2+
  \eps_2\norm{f(x)^{2\over3}u}_{L^2}^2+C_{\eps_1,\eps_2}\norm{u}_{L^2}^2.
\end{eqnarray}
To prove \reff{reqA1}, we use \reff{cf+}; this gives
\[
\abs{\com{P,~f(x)^{1\over3}}}\leq C \abs{y}f(x)^{s-{2\over3}},
\]
and hence
\[
 \forall~\eps_1>0, \quad{\rm Re}\comi{\com{P,~f(x)^{1\over3}}u,~f(x)^{1\over3}u}_{L^2}\leq
 \eps_1\norm{\Lambda_{y}f(x)^{1\over3}u}_{L^2}^2+C_{\eps_1}\norm{f(x)^{s-{2\over3}}u}_{L^2}^2.
\]
Since $s-{2\over3}<{2\over3}$ for $s<{4\over3}$ then the following
interpolation inequality holds:
\[
  \forall~\eps_2>0, \quad \norm{f(x)^{s-{2\over3}}u}_{L^2}^2\leq
  \eps_2\norm{f(x)^{2\over3}u}_{L^2}^2+C_{\eps_2}\norm{u}_{L^2}^2.
\]
Now combination of the above inequalities yields \reff{reqA1}.

{\bf \emph{Step II}.}~Next we will show that  there exists a
constant $C>0$ such that
\begin{eqnarray}\label{reqB}
  \forall~u\in C_0^\infty\big(\mathbb{R}^{2n}\big),\quad \norm{f(x)^{{2\over3}}
  u}_{L^2}^2\leq C\set{{\rm Re}\comi{iQ u,~Ru}_{L^2}+\norm{Pu}_{L^2}^2+\norm{u}_{L^2}^2}.
\end{eqnarray}
Since $Q=y \cdot D_x-\partial_x V(x)\cdot
 D_y$ and $R=2f(x)^{-{2\over3}}\partial_x V(x)\cdot y$, then it's a straightforward verification to see that
 \begin{eqnarray*}
  \frac{i}{2}\Big[R,~Q\Big]=f(x)^{-{2\over3}}\abs{\partial_x
  V(x)}^2-y\cdot\partial_x\inner{f(x)^{-{2\over3}}\partial_x V(x)\cdot
  y}.
 \end{eqnarray*}
 As a result, we use the relation $${\rm Re}\comi{iQ u,~Ru}_{L^2}=
  \frac{i}{2}\comi{\com{R,~Q_{x_0}}u,~u}_{L^2}$$ to get
\begin{eqnarray*}
  {\rm Re}\comi{iQ u,~Ru}_{L^2}=\norm{f(x)^{2\over3}u}_{L^2}^2-
  \norm{f(x)^{-{1\over3}}u}_{L^2}^2-\comi{\inner{y\cdot\partial_x\Big(f(x)^{-{2\over3}}\partial_x V(x)\cdot
  y\Big)}u,~u}_{L^2}.
\end{eqnarray*}
This gives
\begin{eqnarray*}
  \norm{f(x)^{{2\over3}}u}_{L^2}^2\leq {\rm Re}\comi{iQ u,~Ru}_{L^2}+
  \norm{u}_{L^2}^2+\comi{\abs{y\cdot\partial_x\inner{f(x)^{-{2\over3}}\partial_x V(x)\cdot
  y}}u,~u}_{L^2}.
\end{eqnarray*}
Moreover, by use of \reff{cf+}, we compute
\[
  \abs{y\cdot\partial_x\inner{f(x)^{-{2\over3}}\partial_x V(x)\cdot
  y}}\leq C f(x)^{s-{2\over3}}\abs{y}^2\leq C f(x)^{{2\over3}}\abs{y}^2,
\]
which implies that for any $\eps>0,$
\begin{eqnarray*}
  \comi{\abs{y\cdot\partial_x\inner{f(x)^{-{2\over3}}\partial_x V(x)\cdot
  y}}u,~u}_{L^2}&\leq& C \norm{\abs{y}f(x)^{{1\over3}}u}_{L^2}^2
  \leq C \norm{\Lambda_{y}f(x)^{{1\over3}}u}_{L^2}^2
  \\&\leq& \eps \norm{f(x)^{{2\over3}}
  u}_{L^2}^2+C_\eps\set{\norm{Pu}_{L^2}^2+\norm{u}_{L^2}^2},
\end{eqnarray*}
the last inequality using \reff{reqA}. Consequently,
\begin{eqnarray*}
  \forall~\eps>0,\quad\norm{f(x)^{{2\over3}}u}_{L^2}^2\leq {\rm Re}\comi{iQ u,~Ru}_{L^2}+\eps \norm{f(x)^{{2\over3}}
  u}_{L^2}^2+C_\eps\set{\norm{Pu}_{L^2}^2+\norm{u}_{L^2}^2}.
\end{eqnarray*}
Letting $\eps>0$ small enough gives \reff{reqB}.

{\bf \emph{Step III}.}~Now we prove that for any $\eps>0$ there
exists a constant $C_\eps$ such that
\begin{eqnarray}\label{reqC}
  \forall~u\in C_0^\infty\big(\mathbb{R}^{2n}\big),\quad
  \sum_{j=1}^n\abs{\comi{L_j^*L_ju,~Ru}_{L^2}}
  \leq \eps \norm{f(x)^{{2\over3}}
  u}_{L^2}^2+C_\eps\set{\norm{Pu}_{L^2}^2+\norm{u}_{L^2}^2}.
\end{eqnarray}
As a preliminary step,  we firstly show the following estimate:
\begin{eqnarray}\label{reqC1}
  \forall~\eps>0,\quad \norm{\comi{y}^2u}_{L^2}^2
  \leq \eps
  \norm{f(x)^{2\over3}u}_{L^2}+C_\eps\set{\norm{Pu}_{L^2}^2+\norm{u}_{L^2}^2},
\end{eqnarray}
where $\comi{y}=\inner{1+\abs{y}^2}^{1\over2}$. Using \reff{velo}
gives
\begin{eqnarray*}
  \norm{\comi{y}^2u}_{L^2}^2
  &\leq& C\set{{\rm
  Re}\comi{P\comi{y}u,~\comi{y}u}_{L^2}+\norm{\comi{y}u}_{L^2}^2}\\
  &=& C\set{{\rm Re}\comi{Pu,~\comi{y}^{2}u}_{L^2}+
  {\rm
  Re}\comi{\com{P,~\comi{y}}u,~\comi{y}u}_{L^2}}+C\norm{\comi{y}u}_{L^2}^2.
\end{eqnarray*}
This together with \reff{velo} implies that
\begin{eqnarray}\label{090930+1}
  \norm{\comi{y}^2u}_{L^2}^2
  \leq C\set{\norm{Pu}_{L^2}^2+\norm{u}_{L^2}^2}+C
  \abs{\comi{\com{P,~\comi{y}}u,~\comi{y}u}_{L^2}}.
\end{eqnarray}
Moreover observe that
\[
\abs{\com{P,~\comi{y}}u}\leq
C\set{\abs{\partial_xV(x)}\abs{u}+\abs{\partial_yu}+\abs{u}}\leq
C\set{f(x)\abs{u}+\abs{\partial_yu}+\abs{u}},
\]
and hence for any $\eps>0,$
\begin{eqnarray*}
  \abs{\comi{\com{P,~\comi{y}}u,~\comi{y}u}_{L^2}}
  \leq \eps \norm{f(x)^{2\over3}u}_{L^2}^2+C_\eps\set{\norm{\Lambda_yf(x)^{1\over3}u}_{L^2}^2+
  \norm{\Lambda_yu}_{L^2}^2}.
\end{eqnarray*}
This along with \reff{reqA} and \reff{velo} gives
\begin{eqnarray}\label{090930+2}
  \forall~\eps>0,\quad\abs{\comi{\com{P,~\comi{y}}u,~\comi{y}u}_{L^2}}\leq \eps
  \norm{f(x)^{2\over3}u}_{L^2}^2+C_\eps\set{\norm{Pu}_{L^2}^2+\norm{u}_{L^2}^2}.
\end{eqnarray}
Now combining \reff{090930+1} and \reff{090930+2}, we get
\reff{reqC1}.  As a result of \reff{reqC1}, we have
\begin{eqnarray}\label{reqC2}
  \forall~\eps>0,\quad \norm{\Lambda_y\comi{y}u}_{L^2}^2
  \leq \eps
  \norm{f(x)^{2\over3}u}_{L^2}^2+C_\eps\set{\norm{Pu}_{L^2}^2+\norm{u}_{L^2}^2}.
\end{eqnarray}
Indeed by \reff{velo} one has
\begin{eqnarray*}
   \norm{\Lambda_y\comi{y}u}_{L^2}^2
  &\leq& C\set{{\rm Re}\comi{P\comi{y}u,~\comi{y}u}_{L^2}+\norm{\comi{y}u}_{L^2}^2}\\
  &\leq&C\abs{\comi{\com{P,~\comi{y}}u,~\comi{y}u}_{L^2}}+C\set{\norm{Pu}_{L^2}^2+\norm{\comi{y}^{2}u}_{L^2}}.
\end{eqnarray*}
So \reff{reqC2} can be deduced from \reff{reqC1} and
\reff{090930+2}. Now we are ready to prove \reff{reqC}. Observe
\begin{eqnarray*}
  \abs{\comi{L_j^*L_ju,~Ru}_{L^2}}=\abs{\comi{f(x)^{{1\over3}}L_ju,~f(x)^{-{1\over3}}L_jRu}_{L^2}}\leq
  \norm{\Lambda_yf(x)^{{1\over3}}u}_{L^2}^2+\norm{f(x)^{-{1\over3}}L_jRu}_{L^2}^2.
\end{eqnarray*}
Then in view of \reff{reqA}, we see that the required inequality
\reff{reqC} will follow if the following estimate holds:
\begin{eqnarray}\label{reqC3}
  \forall~\eps>0,\quad \norm{f(x)^{-{1\over3}}L_jRu}_{L^2}^2\leq\eps
  \norm{f(x)^{2\over3}u}_{L^2}^2+C_\eps\set{\norm{Pu}_{L^2}^2+\norm{u}_{L^2}^2}.
\end{eqnarray}
Since
\[
  L_jRu=2uf(x)^{-{2\over3}}\partial_{y_j}\inner{\partial_x V(x)\cdot
  y}+2f(x)^{-{2\over3}}\inner{\partial_x V(x)\cdot
  y}\partial_{y_j}u+f(x)^{-{2\over3}}y_j\inner{\partial_x V(x)\cdot
  y}u,
\]
then
\begin{eqnarray*}
  \norm{f(x)^{-{1\over3}}L_jRu}_{L^2}^2\leq C\set{\norm{u}_{L^2}^2
  +\norm{\Lambda_y\comi{y}u}_{L^2}^2}.
\end{eqnarray*}
This along with  \reff{reqC2} gives \reff{reqC3}, completing the
proof \reff{reqC}.

Now we combine the inequalities \reff{rmain}, \reff{reqA},
\reff{reqB} and \reff{reqC}, to obtain
\begin{eqnarray*}
  \forall~\eps>0,\quad \norm{f(x)^{{2\over3}}
  u}_{L^2}^2\leq \eps \norm{f(x)^{{2\over3}}
  u}_{L^2}^2+C_\eps\set{\norm{Pu}_{L^2}^2+\norm{u}_{L^2}^2}.
\end{eqnarray*}
Taking $\eps=\frac{1}{2}$ gives the desired estimate \reff{req}.
This completes the proof of Proposition \ref{prop+1}. \qed

\section{Proof of  Theorem \ref{+Hypo}: the second part}
\label{sec4}

This section is devoted to the proof of \reff{+++A1}, and then the
proof of Theorem \ref{+Hypo} will be completed. As a convention, we
use the capital letter $C$ to denote different suitable constants.
Let $V$ satisfy the assumption \reff{+HyP}. In the sequel we use the
notation
\[f(x)=\inner{1+\abs{\partial_xV(x)}^2}^{1\over2}.\]
Then  \reff{+HyP} yields
   \begin{eqnarray}\label{091005+1}
     \forall~x\in\mathbb{R}^{n},\quad \abs{\partial_xf(x)}\leq C f(x)^s.
   \end{eqnarray}
In view of \reff{velo}, to prove \reff{+++A1} one only has to show

\begin{proposition}\label{prpfin}
If $V(x)$ satisfies the assumption \reff{+HyP}, then
\begin{eqnarray}\label{++A1}
   \forall~u\in C_0^\infty\big(\mathbb{R}^{2n}\big),\quad \norm{\inner{1-\triangle_x}^{\frac{\delta}{2}}u}_{L^2}
   \leq C\Big\{
   \norm{P u}_{L^2}+\norm{u}_{L^2}\Big\},
\end{eqnarray}
where $\delta$  equals to $\frac{2}{3}$ if $s\leq \frac{2}{3}$,
$\frac{4}{3}-s$ if ${2\over3}<s\leq \frac{10}{9},$ and
$\frac{2}{3}-{s\over2}$ if ${{10}\over9}< s<\frac{4}{3}.$
\end{proposition}

We will use localization arguments to prove the above proposition.
Firstly let's  recall some standard results concerning the partition
of unity. For more detail we refer to \cite{Hormander85} for
instant. Let $g$ be a metric of the following form
\begin{eqnarray}\label{7292}
  g_x=f(x)^{s}\abs{dx}^2,\quad x\in\mathbb{R}^n,
\end{eqnarray}
where $s$ is  the real number given in \reff{091005+1}.

\begin{lemma}\label{lem729}
  Suppose $f$ satisfies the assumption \reff{091005+1}. Then the metric
  $g$ defined by \reff{7292} is slowly varying,  i.e., we can find two constants $C_*,r>0$ such
  that if $g_x(x-y)\leq r^2$ then $$C_*^{-2}\leq \frac{g_x}{g_y}\leq
  C_*^2.$$
\end{lemma}

\noindent\emph{Proof.}
  We only need to show that
  \begin{eqnarray}\label{slow+}
    \exists~ r,C_*>0,~\forall~x,y\in\mathbb{R}^n,~~\quad \abs{x-y}\leq r f(x)^{-{s\over2}}\Longrightarrow
    C_*^{-1} \leq \frac{f(x)^{s\over2}}{f(y)^{s\over2}}\leq C_*.
  \end{eqnarray}
  Making use of \reff{cf+} and the fact that $s<\frac{4}{3},$ we
  have
  \begin{eqnarray*}
    \forall~x\in\mathbb{R}^n,\quad\abs{\partial_{x}\inner{f(x)^{-{s\over2}}}}\leq
    f(x)^{-{s\over2}-1}\abs{\partial_{x}f(x)}\leq
    C f(x)^{{s\over2}-1}\leq  C
  \end{eqnarray*}
  with $C$ the constant in \reff{cf+}. As a consequence, one can
  find a
  constant $\tilde C$ depending only on $C$ and the dimension $n,$
  such that
  \begin{eqnarray*}
     \forall~x,y\in\mathbb{R}^n, \quad\abs{f(x)^{-{s\over2}}-f(y)^{-{s\over2}}}\leq \tilde C\abs{x-y},
  \end{eqnarray*}
  from which we conclude that if $\abs{x-y}\leq r f(x)^{-{s\over2}}$ then
  \begin{eqnarray*}
     \abs{f(x)^{-{s\over2}}-f(y)^{-{s\over2}}}\leq  r\tilde C f(x)^{-{s\over2}}.
  \end{eqnarray*}
  Thus
  \begin{eqnarray*}
     \abs{\frac{f(x)^{{s\over2}}}{f(y)^{{s\over2}}}-1}\leq  r\tilde C.
  \end{eqnarray*}
  This gives   \reff{slow+} if we choose $r=\frac{\tilde C}{2}$ and $C_*=2.$
\qed

Let $g$ be the   metric given by \reff{7292}. We denote by $S(1,g)$
the class of smooth real-valued functions $a(x)$ satisfying
  the following condition:
  \begin{eqnarray*}
    \forall~\gamma\in \mathbb{Z}_+^n,~~\forall~x\in\mathbb{R}^n,\quad
    \abs{\partial^\gamma a(x)}\leq C_\gamma f(x)^{\frac{s\abs\gamma}{2}}.
  \end{eqnarray*}
  The space $S(1,g)$ endowed with the seminorms
  \[
    \abs{a}_{k,S(1,g)}=\sup_{x\in\mathbb{R}^n,\abs{\gamma}=k}f(x)^{-\frac{sk}{2}}\abs{\partial^\gamma a(x)},\quad
    k\geq 0,
  \]
  becomes a Fr\'{e}chet space.

The main feature of a slowly varying metric is that it allows us to
introduce some partitions of unity related to the metric.  We state
it as the following lemma.

\begin{lemma}[(Lemma 18.4.4. of \cite{Hormander85})]\label{uni}
  Let $g$ be a slowly varying metric. We can find a constant $r_0>0$ and a sequence $x_\mu\in\mathbb{R}^n,\mu\geq1,$ such
  that the union of the balls
  \[
    \Omega_{\mu,r_0}=\set{x\in\mathbb{R}^{n};\quad g_{x_\mu}\inner{x-x_\mu}<r_0^2}
  \]
  coves the whole space $\mathbb{R}^{n}.$ Moreover there exists a positive integer
  $N ,$ depending only on $r_0,$ such that the intersection
  of more than $N $ balls is always empty.
  One can choose a
  family of nonnegative functions $\set{\varphi_\mu}_{\mu\geq 1}$ uniformly bounded in
  $S(1,g)$ such that
 \begin{eqnarray}\label{73966}
   {\rm supp}~ \varphi_\mu\subset\Omega_{\mu,r_0},\quad
   \sum_{\mu\geq1} \varphi_\mu^2 =1~~\,\,{\rm and }~~\,\,
   \sup_{\mu\geq1}\abs{\partial_x \varphi_\mu (x)}\leq C f(x)^{s\over2}.
\end{eqnarray}
  Here by  uniformly bounded in
  $S(1,g),$ we mean
  \begin{eqnarray*}
     \sup_{\mu}\abs{\varphi_\mu}_{k,S(1,g)}\leq C_{k}, \quad k\geq
     0.
  \end{eqnarray*}
\end{lemma}

\begin{remark}\label{7308}
  If we choose $r_0$ small enough such that $r_0\leq r$ with $r$ the constant given in Lemma \ref{lem729},
  then there exists a constant $C,$ such that for any
  $\mu\geq1$ one has
\begin{eqnarray}\label{7310}
  \forall~x,y\in {\rm
  supp}~  \varphi_\mu,\quad C^{-1}f(y) \leq f(x) \leq C f(y).
\end{eqnarray}
\end{remark}

\begin{lemma}
  Let $V(x)$ satisfy the assumption \reff{+HyP},  and
  let $\set{\varphi_\mu}_{\mu\geq1}$ be the partition of unity given
  above. Then we have for any $u\in  C_0^\infty\big(\mathbb{R}^{2n}\big),$
  \begin{eqnarray}\label{AB1}
    \sum\limits_{\mu\geq1}\norm{\inner{y\cdot\partial_x\varphi_\mu}
    u}_{L^2}^2\leq C
    \norm{ \Lambda_yf(x)^{s\over2}u}_{L^2}^2
  \end{eqnarray}
  and
  \begin{eqnarray}\label{AB2}
    \sum\limits_{\mu\geq1}\norm{\varphi_\mu(x)\inner{\partial_x V(x)-\partial_x
    V(x_\mu)}\cdot \partial_y
    u}_{L^2}^2 \leq C
    \norm{\Lambda_y f(x)^{s\over2}u}_{L^2}^2.
  \end{eqnarray}
\end{lemma}

\noindent\emph{Proof.} Firstly we show \reff{AB1}. Observe
\begin{eqnarray*}
  \norm{\inner{y\cdot\partial_x\varphi_\mu}u}_{L^2}^2
  =\comi{\inner{y\cdot\partial_x\varphi_\mu}^2u,~u}_{L^2},
\end{eqnarray*}
and by Lemma \ref{uni}, we see that  $\sum_{\mu\geq1}
\abs{\partial_x\varphi_\mu}^2$ is a sum of at most $N$ terms and
hence bounded from above by $f^{s}.$ As a result,
  \[
    \sum_{\mu\geq1} \inner{y\cdot\partial_x
    \varphi_\mu(x)}^2\leq C \abs{y}^2 \sum_{\mu\geq1} \abs{\partial_x\varphi_\mu}^2
    \leq C\abs{y}^2 f^{s} .
  \]
  Then \reff{AB1} follows. Next we estimate \reff{AB2}.
 Note that $\abs{x-x_\mu}\leq Cf(x_\mu)^{-{s\over2}}$ for any $x\in$
 supp~$\varphi_\mu,$
   and hence we can deduce from \reff{+HyP} and \reff{7310}
    that
   \[
     \sum_{\mu\geq 1}\varphi_\mu(x)^2\abs{\partial_x V(x)-\partial_x
     V(x_\mu)}^2\leq C\sum_{\mu\geq 1}\varphi_\mu(x)^2
     f(x)^{2s}\abs{x-x_\mu}
     \leq C\sum_{\mu\geq 1}\varphi_\mu(x)^2 f(x)^{s}\leq C f(x)^{s}.
   \]
   This along with the inequality
   \[
   \sum\limits_{\mu\geq1}\norm{\varphi_\mu(x)\inner{\partial_x V(x)-\partial_x
    V(x_\mu)}\cdot \partial_y
    u}_{L^2}^2=\comi{\sum_{\mu\geq 1}\varphi_\mu(x)^2\abs{\partial_x V(x)-\partial_x
     V(x_\mu)}^2\abs{\partial_y u},\abs{\partial_y u}}_{L^2}
   \]
    implies \reff{AB2}. Then the proof is completed.
  \qed

  \begin{lemma}\label{lemreg}
     Let $\set{\varphi_\mu}_{\mu\geq1}$ be the partition given in Lemma \ref{uni},
     and let $a\in]0, 1/2[$ be a real number.
     Then there exists a constant $C,$ depending on the integer $N $
     given in Lemma \ref{uni}, such that
     \begin{equation}\label{tra}
       \forall~u\in C_0^\infty\big(\mathbb{R}^{2n}\big),\quad
       \norm{\inner{1-\triangle_x}^{a}u}_{L^2}^2\leq
       C\sum_{\mu\geq1}\norm{\inner{1-\triangle_x}^{a}\varphi_\mu  u}_{L^2}^2
       +C\norm{Pu}_{L^2}^2+C\norm{u}_{L^2}^2.
     \end{equation}
   \end{lemma}

   In order to prove Lemma \ref{lemreg} we need  the following technical lemma.

\begin{lemma}
  Let $b\in]0,1[$ be a real number and $\abs{D_x}^b$ be the Fourier multiplier defined
  by, with $u\in C_0^\infty(\mathbb R^n),$
  \[
    \abs{D_x}^b u(x)=\mathcal {F}^{-1}\inner{\abs \xi^b \hat
    u(\xi)}.
  \]
  Let $\set{\varphi_\mu}_{\mu\geq1}$ be the partition given in Lemma
  \ref{uni}.
  Then there exists a constant $C$  such that for any $u\in C_0^\infty(\mathbb R^n),$
  \begin{eqnarray}\label{2010090302}
    \norm{\sum_{\mu\geq1}\com{|D_x|^b,\:f^{-s/2}\varphi_\mu}\varphi_\mu u}_{L^2(\mathbb R^n)} \leq C\norm{u}_{L^2(\mathbb
    R^n)}
  \end{eqnarray}
  and
  \begin{eqnarray}\label{2010090701}
    \norm{ \com{|D_x|^b,\:f^{-s/2}}  u}_{L^2(\mathbb R^n)} \leq C\norm{u}_{L^2(\mathbb
    R^n)}.
  \end{eqnarray}
  Recall here $f(x)=\inner{1+\abs{\partial_x V(x)}^2}^{1/2}$ and $s$
  is the real number given in \reff{091005+1}.
\end{lemma}

\noindent\emph{Proof.}  In the proof we use $C$ to denote different
suitable positive constants, and for simplicity we use the notation
\[
   \omega_\mu=f^{-s/2}\varphi_\mu.
\]
In view of Lemma \ref{uni} and the estimate \reff{091005+1}, we have
\begin{equation}\label{20100903}
    \sup_{x\in\mathbb R^n}\inner{\sum_{\mu\geq 1}\abs{\varphi_\mu(x)}^2}^{1\over2}+
    \sup_{x,x'\in\mathbb R^n}\inner{\sum_{\mu\geq 1}\abs{\omega_\mu(x)-\omega_\mu(x')}^2}^{1\over2}+\sup_{x\in\mathbb
    R^n}\inner{\sum_{\mu\geq1}\abs{\partial_x \omega_\mu (x)}^2}^{1\over2}\leq C.
\end{equation}
Next we will show the following relation
\begin{eqnarray}\label{FourierM}
  \forall~u\in C_0^\infty(\mathbb{R}^n),\quad |D_x|^b\, u(x)=C_b\int_{\mathbb R^n}\frac{u(x)-u(x-\tilde
 x)}{\abs
  {\tilde x}^{n+\sigma}}\,d\tilde x
\end{eqnarray}
with $C_b\neq 0$ being a complex constant depending only on the real
number $b$ and the dimension $n.$ In fact, the inverse Fourier
transform implies
\[
  \int_{\mathbb R^n}\frac{u(x)-u(x-\tilde x)}{\abs
  {\tilde x}^{n+b}}\,d\tilde x=\int_{\mathbb R^n} \hat u(\xi)\,e^{i\,x\cdot\xi}\left(
  \int_{\mathbb R^n}\frac{1-e^{-i\,\tilde x\cdot\xi}}{\abs
  {\tilde x}^{n+b}}\,d\tilde x\right)\,d\xi
\]
On the other hand, we can verify that
\[
  \int_{\mathbb R^n}\frac{1-e^{-i\,\tilde x\cdot\xi}}{\abs
  {\tilde x}^{n+b}}\,d\tilde x=\abs\xi^b \int_{\mathbb R^n}
  \frac{1-e^{-i\, z\cdot\frac{\xi}{\abs\xi}}}{\abs
  {z}^{n+b}}\,d z.
\]
Observe that $\int_{\mathbb R^n}
  \frac{1-e^{-i\, z\cdot\frac{\xi}{\abs\xi}}}{\abs
  {z}^{n+b}}\,d z \neq 0$ is a complex constant depending only on
 $b$ and
the dimension $n$, but independent of $\xi.$  Then the above two
equalities give \reff{FourierM}. Now we use \reff{FourierM} to get
\begin{equation*}
\begin{split}
  |D_x|^b\,\biginner{\omega_\mu\,\varphi_\mu\, u}(x)&=
  C_b\int_{\mathbb R^n}\frac{\omega_\mu(x)\varphi_\mu(x)u(x)
  -\omega_\mu(x-\tilde x)\varphi_\mu(x-\tilde x)u(x-\tilde x)}{\abs
  {\tilde x}^{n+b}}\,d\tilde x\\&= \omega_\mu(x)|D_x|^b\, (\varphi_\mu\, u)(x)
  +C_b\int_{\mathbb R^n}\frac{
  \varphi_\mu(x-\tilde x)u(x-\tilde
  x)\biginner{\omega_\mu(x)
  -\omega_\mu(x-\tilde x)}}{\abs
  {\tilde x}^{n+b}}\,d\tilde x,
\end{split}
\end{equation*}
which gives
\begin{eqnarray}\label{formula}
  \com{|D_x|^b,\:\omega_\mu}(\varphi_\mu\, u)(x)=
  C_b\int_{\mathbb R^n}\frac{\varphi_\mu(x-\tilde x)u(x-\tilde
  x)\biginner{\omega_\mu(x)
  -\omega_\mu(x-\tilde x)}}{\abs
  {\tilde x}^{n+b}}\,d\tilde x.
\end{eqnarray}
Let $\rho$ be the characteristic function of the unit ball
$\set{x\in\mathbb R^n;\: \abs x\leq 1}.$ We compute
\begin{equation*}
\begin{split}
   &\norm{\sum_{\mu\geq1}\com{|D_x|^b,\,\omega_\mu}\varphi_\mu u}_{L^2}
   ^2=
   |C_b| ^2\int_{\mathbb R^n}\inner{\sum_{\mu\geq1}\int_{\mathbb R^n}\frac{u(x-\tilde
 x)\varphi_\mu(x-\tilde x)\biginner{\omega_\mu(x)
  -\omega_\mu(x-\tilde x)}}{\abs
  {\tilde x}^{n+b}}d\tilde x}^2dx\\
   &\leq
 2 |C_b| ^2\int_{\mathbb R^n}\inner{\sum_{\mu\geq1}\int_{\mathbb R^n}\frac{ \rho(\tilde x)
   u(x-\tilde x)\varphi_\mu(x-\tilde x)\biginner{\omega_\mu(x)
  -\omega_\mu(x-\tilde x)}}{\abs
  {\tilde x}^{n+b}}\,d\tilde x}^2dx\\
   &\quad +2 |C_b| ^2\int_{\mathbb R^n}\inner{\sum_{\mu\geq1}\int_{\mathbb R^n}\frac{\inner{1-
   \rho (\tilde x)}
   u(x-\tilde x)\varphi_\mu(x-\tilde x)\biginner{\omega_\mu(x)
  -\omega_\mu(x-\tilde x)}}{\abs
  {\tilde x}^{n+b}}\,d\tilde x}^2dx\\
  &=: \mathcal A_1+\mathcal A_2.
\end{split}
\end{equation*}
Now we treat the terms $\mathcal A_1$ and $\mathcal A_2.$  Cauchy's
inequality yields
\begin{equation*}
\abs{\sum_{\mu\geq 1} \varphi_\mu(x-\tilde
x)\inner{\omega_\mu(x)-\omega_\mu(x-\tilde x)}}\leq
\inner{\sum_{\mu\geq 1}\abs{\varphi_\mu(x-\tilde
x)}^2}^{1\over2}\inner{\sum_{\mu\geq
1}\abs{\omega_\mu(x)-\omega_\mu(x-\tilde x)}^2}^{1\over2}.
\end{equation*}
This along with \reff{20100903} gives that for any $x,\tilde
x\in\mathbb R^n,$  we have $\abs{\sum_{\mu\geq 1}
\varphi_\mu(x-\tilde x)\inner{\omega_\mu(x)-\omega_\mu(x-\tilde
x)}}\leq C $ and hence
\[
   \mathcal {A}_2\leq C
\int_{\mathbb R^n}\inner{\int_{\mathbb R^n}\frac{\inner{1-
\rho(\tilde x)}
  \abs{u(x-\tilde x)}}{\abs
  {\tilde x}^{n+b}}\,d\tilde x}^2dx.
\]
Moreover, using the relation
\[
 \omega_\mu(x)-\omega_\mu(x-\tilde x)= \int_0^1 \partial_{x}
 \omega_{\mu}\biginner{tx+(1-t)(x-\tilde x)}\cdot\tilde x ~dt
\]
and the inequality \reff{20100903} yields that for any $x,\tilde
x\in\mathbb R^n$  we have $\inner{\sum_{\mu\geq
1}\abs{\omega_\mu(x)-\omega_\mu(x-\tilde x)}^2}^{1\over2}\leq C
\abs{\tilde x}$ and hence
\[
\abs{\sum_{\mu\geq 1} \varphi_\mu(x-\tilde
x)\inner{\omega_\mu(x)-\omega_\mu(x-\tilde x)}} \leq C \abs{\tilde
x},
\]
which implies
\[
   \mathcal {A}_1\leq C
\int_{\mathbb R^n}\inner{\int_{\mathbb R^n}\frac{  \rho(\tilde x)
  \abs{u(x-\tilde x)}}{\abs
  {\tilde x}^{n+b-1}}\,d\tilde x}^2dx.
\]
Combining these inequalities gives
\begin{equation*}
\begin{split}
   \norm{\sum_{\mu\geq1}\com{|D_x|^b,\:\omega_\mu}\varphi_\mu u}_{L^2}^2&\leq  C
\int_{\mathbb R^n}\inner{\int_{\mathbb R^n}\frac{ \rho(\tilde x)
  \abs{u(x-\tilde x)}}{\abs
  {\tilde x}^{n+b-1}}\,d\tilde x}^2dx\\
  &\quad+C
\int_{\mathbb R^n}\inner{\int_{\mathbb R^n}\frac{\inner{1-
\rho(\tilde x)}
  \abs{u(x-\tilde x)}}{\abs
  {\tilde x}^{n+b}}\,d\tilde x}^2dx.
\end{split}
\end{equation*}
Moreover, for the terms on the right side of the above inequality,
we can use Young's inequality for convolutions and the fact that
$\rho$ is the characteristic function of the unit ball, to get
\[
\int_{\mathbb R^n}\inner{\int_{\mathbb R^n}\frac{ \rho(\tilde x)
  \abs{u(x-\tilde x)}}{\abs
  {\tilde x}^{n+b-1}}\,d\tilde x}^2dx\leq C
  \norm{u}_{L^2(\mathbb R^n)}^2\Big\|\frac{ \rho }{\abs
  { x}^{n+b-1}}\Big\|_{L^1(\mathbb R^n) }^2\leq
   C\norm{u}_{L^2(\mathbb R^n)}^2
\]
and
\[
\int_{\mathbb R^n}\inner{\int_{\mathbb R^n}\frac{\inner{1-
\rho(\tilde x)}
  \abs{u(x-\tilde x)}}{\abs
  {\tilde x}^{n+b}}\,d\tilde x}^2dx\leq
  C\norm{u}_{L^2(\mathbb R^n)}^2\Big\|\frac{ 1-\rho }{\abs
  { x}^{n+b}}\Big\|_{L^1(\mathbb R^n) }^2\leq   C
  \norm{u}_{L^2(\mathbb R^n)}^2.
\]
We combine these inequalities to  get the desired estimate
\reff{2010090302}. The estimate \reff{2010090701}, which is easier
to treat, can be obtained via the similar arguments as above. This
completes the proof. \qed

\bigskip
   \noindent\emph{Proof of Lemma \ref{lemreg}.} We only need show
   that, with $b\in]0,1[,$
   \begin{eqnarray}\label{tra++}
       \forall~u\in C_0^\infty\big(\mathbb{R}^{2n}\big),\quad
       \norm{\abs{D_x}^{b}u}_{L^2}^2\leq
       C\sum_{\mu\geq1}\norm{\abs{D_x}^{b}\varphi_\mu  u}_{L^2}^2
       +C\norm{Pu}_{L^2}^2+C\norm{u}_{L^2}^2.
     \end{eqnarray}
     By \reff{73966}, we see
     $\norm{\abs{D_x}^{b}u}_{L^2}^2=\norm{\sum_{\mu\geq 1}\abs{D_x}^{b}\varphi_\mu^2\,u}_{L^2}^2.$ Thus
     \begin{eqnarray}\label{100801}
       \norm{\abs{D_x}^{b}u}_{L^2}^2
       \leq2\norm{\sum_{\mu\geq 1}\com{\abs{D_x}^{b},~f^{-{s\over2}}\varphi_\mu}\varphi_\mu f^{{s\over2}}u}_{L^2}^2+
       2\norm{\sum_{\mu\geq
       1}f^{-{s\over2}}\varphi_\mu\abs{D_x}^{b}\varphi_\mu\,f^{{s\over2}}u}_{L^2}^2.
     \end{eqnarray}
     In view of \reff{2010090302} we have
     \begin{eqnarray}\label{100802}
       \norm{\sum_{\mu\geq 1}\com{\abs{D_x}^{b},~f^{-{s\over2}}\varphi_\mu}\varphi_\mu f^{{s\over2}}u}_{L^2}^2
       \leq C\norm{f^{{s\over2}}u}_{L^2}^2 \leq
       C\norm{P\,u}_{L^2}^2+C\norm{u}_{L^2}^2,
     \end{eqnarray}
     the last inequality following from \reff{req}. It remains to
     handle the second term on the right side of \reff{100801}.
     For each $\mu\geq1,$ set
     \[
       I_\mu=\set{\nu\geq 1;\quad {\rm supp}~\varphi_{\nu}\cap~{\rm
       supp}~\varphi_{\mu}\neq\emptyset}.
     \]
     Then $I_\mu$ is a finite set and has at most $N $ elements.  Recall $N $ is the integer
     given in Lemma \ref{uni} such that
     the intersection of more than $N $ balls is always empty.
     Direct calculus give that for any $u\in C_0^\infty(\mathbb{R}^{2n}),$
     \begin{equation*}
     \begin{split}
       &\norm{\sum_{\mu\geq
       1}f^{-{s\over2}}\varphi_\mu\abs{D_x}^{b}\varphi_\mu\,f^{{s\over2}}u}_{L^2}^2=\sum_{\mu\geq1}
       \sum_{\nu\in I_\mu}\comi{\varphi_\mu f^{-{s\over2}}
      \abs{D_x}^b \varphi_\mu f^{{s\over2}}u,~~\varphi_\nu f^{-{s\over2}}
      \abs{D_x}^b \varphi_\nu f^{{s\over2}}u}_{L^2}\\
      &\leq \sum_{\mu\geq1}\sum_{\nu\in I_\mu}\norm{\varphi_\mu f^{-{s\over2}}
      \abs{D_x}^b \varphi_\mu f^{{s\over2}}u}_{L^2}^2+\sum_{\mu\geq1}\sum_{\nu\in I_\mu}
      \norm{\varphi_\nu f^{-{s\over2}}
      \abs{D_x}^b \varphi_\nu f^{{s\over2}}u}_{L^2}^2\\
      &=2\sum_{\mu\geq1}\sum_{\nu\in I_\mu}\norm{\varphi_\mu f^{-{s\over2}}
      \abs{D_x}^b \varphi_\mu f^{{s\over2}}u}_{L^2}^2\\
      \\
      &\leq2\sum_{\mu\geq1}\sum_{\nu\in I_\mu}\norm{ f^{-{s\over2}}
      \abs{D_x}^b \varphi_\mu f^{{s\over2}}u}_{L^2}^2.
    \end{split}
    \end{equation*}
     Since $I_\mu$ has at most $N$ elements then  it follows that
     \begin{eqnarray}\label{tra+2}
      \forall~u\in C_0^\infty(\mathbb{R}^{2n}), \quad
      \norm{\sum_{\mu\geq
       1}f^{-{s\over2}}\varphi_\mu\abs{D_x}^{b}\varphi_\mu\,f^{{s\over2}}u}_{L^2}^2
       \leq 2N \sum_{\mu\geq1}\norm{ f^{-{s\over2}}
      \abs{D_x}^b \varphi_\mu f^{{s\over2}}u}_{L^2}^2.
     \end{eqnarray}
     One the other hand, one can verify
     that
     \begin{eqnarray*}
     \sum_{\mu\geq1}\norm{ f^{-{s\over2}}
      \abs{D_x}^b \varphi_\mu f^{{s\over2}}u}_{L^2}^2&\leq&\sum_{\mu\geq1}\norm{ \com{\abs{D_x}^b,~f^{-{s\over2}}
      } \varphi_\mu f^{{s\over2}}u}_{L^2}^2+\sum_{\mu\geq1}\norm{
      \abs{D_x}^b f^{-{s\over2}}\varphi_\mu f^{{s\over2}}u}_{L^2}^2\\
      &\leq&C\sum_{\mu\geq1}\norm{\varphi_\mu f^{{s\over2}}u}_{L^2}^2+C\sum_{\mu\geq1}\norm{
      \abs{D_x}^b  \varphi_\mu  u}_{L^2}^2\\
      &\leq&C\norm{P\,u}_{L^2}^2+\norm{u}_{L^2}^2+C\sum_{\mu\geq1}\norm{
      \abs{D_x}^b  \varphi_\mu  u}_{L^2}^2,
     \end{eqnarray*}
     the second inequality using \reff{2010090701} and the last inequality using \reff{req}. These inequalities
     along with \reff{tra+2} gives
     \begin{eqnarray*}
      \forall~u\in C_0^\infty(\mathbb{R}^{2n}), \quad
      \norm{\sum_{\mu\geq
       1}f^{-{s\over2}}\varphi_\mu\abs{D_x}^{b}\varphi_\mu\,f^{{s\over2}}u}_{L^2}^2
       \leq C\sum_{\mu\geq1}\norm{
      \abs{D_x}^b  \varphi_\mu  u}_{L^2}^2+C\norm{P\,u}_{L^2}^2+\norm{u}_{L^2}^2 .
     \end{eqnarray*}
     This along with \reff{100801} and \reff{100802} yields the
     desired estimate \reff{tra++} , completing the proof of Lemma \ref{lemreg}.
   \qed

\subsection{End of the proof of Theorem \ref{+Hypo}}\label{sec43}
In this subsection we prove Proposition \ref{prpfin}.  Let
$\set{\varphi_\mu}_{\mu\geq1}$ be the partition of unity given in
Lemma \ref{uni}. For each $\mu\geq1,$  define the operator  $R_\mu$
by
\begin{eqnarray}\label{+Rem}
  R_\mu=-y\cdot\partial_x\varphi_\mu(x)
  -\varphi_\mu\inner{\partial_x V(x)-\partial_x V(x_\mu)}\cdot
  \partial_y.
\end{eqnarray}
We associate with each  $x_\mu\in\mathbb{R}^{n}$  the operator
  \[
    P_{x_\mu}=y\cdot\partial_x-\partial_x V(x_\mu)\cdot
    \partial_y-\triangle_y+\frac{\abs y^2}{4}-\frac{n}{2}.
  \]
  Then  we  have
  \begin{eqnarray*}\label{expression}
     \varphi_\mu P u=P_{x_\mu}~\varphi_\mu\, u+ R_\mu u
  \end{eqnarray*}
  with $R_\mu$ the operator given in  \reff{+Rem}. This gives
  \begin{eqnarray}\label{expression}
     \sum_{\mu\geq 1}\norm{P_{x_\mu} ~\varphi_\mu\,u}_{L^2}^2\leq
    2\sum_{\mu\geq 1}\norm{\varphi_\mu  Pu}_{L^2}^2
    + 2\sum_{\mu\geq 1}\norm{R_\mu u}_{L^2}^2\leq
     2\norm{ Pu}_{L^2}^2
    + 2\sum_{\mu\geq 1}\norm{R_\mu u}_{L^2}^2.
  \end{eqnarray}

  \begin{proposition}\label{PA}
There is a constant $C$ independent of $x_\mu,$  such that
 for any $u\in  C_0^\infty\big(\mathbb{R}^{2n}\big),$ one has
  \begin{eqnarray}\label{+ES}
   \abs{\partial_x
   V(x_\mu)}^{4\over3}\norm{ u}_{L^2}^2+\norm{\inner{1-\triangle_x}^{1\over3} u}_{L^2}^2
   \leq
   C\set{
   \norm{P_{x_\mu} u}_{L^2}^2+\norm{u}_{L^2}^2},
\end{eqnarray}
or equivalently,
\begin{eqnarray}\label{ES}
   \norm{\tilde \Lambda_{x_\mu}^{2\over3}u}_{L^2}^2 \leq
   C\set{
   \norm{P_{x_\mu} u}_{L^2}^2+\norm{u}_{L^2}^2},
\end{eqnarray}
where $\tilde \Lambda_{x_\mu}=\inner{1+\frac{1}{2}\abs{\partial_x
   V(x_\mu)}^2-\triangle_x}^{1\over2}.$

\end{proposition}

The above proposition can be proven in the same way as Proposition
5.22 of \cite{HelfferNier05}, by taking Fourier analysis in the
x-variable and then reducing the problem to a semi-class problem. We
refer to \cite{HelfferNier05} and references therein for more
details. For the sake of completeness we present a direct proof in
the Appendix at the end of the paper.

\medskip
\begin{lemma}\label{patch}
Suppose $V(x)$ satisfies the assumption \reff{+HyP}. Let $R_\mu$ be
the operator given in \reff{+Rem}. Then
  \begin{eqnarray}\label{CMU}
    \forall~u\in   C_0^\infty\big(\mathbb{R}^{2n}\big),\quad\sum\limits_{\mu\geq1}\norm{
    R_\mu
    u}_{L^2}^2 \leq C
    \set{\norm{P f(x)^{\tilde s}
    u}_{L^2}^2+\norm{P
    u}_{L^2}^2+\norm{u}_{L^2}^2},
  \end{eqnarray}
  where $\tilde s=\frac{2}{3}-\delta$ with $\delta$ given in \reff{+++A1}, i.e.,
  $\tilde s$ equals to $0$ if $s\leq {2\over3}$,  $s-{2\over3}$ if ${2\over3}< s\leq\frac{10}{9},$ and
   $\frac{s}{2}$ if ${10\over9}< s<\frac{4}{3}$.
\end{lemma}

\noindent\emph{Proof.}
   As a convention, we use the capital letter $C$ to denote different suitable constants.
   Since $V(x)$ satisfies \reff{+HyP}, then \reff{091005+1} holds. Observe $\sum\limits_{\mu\geq1}\norm{
    R_\mu u}_{L^2}^2$ is bounded from above by
    \[
    2\sum\limits_{\mu\geq1}\norm{\inner{y\cdot\partial_x\varphi_\mu}
    u}_{L^2}^2+2\sum\limits_{\mu\geq1}\norm{\varphi_\mu(x)\inner{\partial_x V(x)-\partial_x
    V(x_\mu)}\cdot \partial_y
    u}_{L^2}^2.
    \]
    Then in view of \reff{AB1} and \reff{AB2}, we have
  \begin{eqnarray*}
    \sum\limits_{\mu\geq1}\norm{R_\mu
    u}_{L^2}^2 \leq C
    \norm{\Lambda_y  f(x)^{s\over2}u}_{L^2}^2.
  \end{eqnarray*}
  So we only have  to treat the term $\norm{\Lambda_y  f(x)^{s\over2}u}_{L^2}^2.$  It follows from \reff{velo}
  that
  \begin{eqnarray*}
    \norm{\Lambda_y  f(x)^{s\over2}u}_{L^2}^2
    \leq C\set{\abs{\comi{
    P f(x)^{s\over2}u,~ f(x)^{s\over2}u}_{L^2}}
    +\norm{ f(x)^{s\over2}u}_{L^2}^2}.
  \end{eqnarray*}
  Since ${s\over2}<\frac{2}{3}$ then by \reff{req} we have
  \begin{eqnarray}\label{0909261}
    \forall~u\in C_0^\infty\inner{\mathbb{R}^{2n}},\quad\norm{ f(x)^{\frac{s}{2}}u}_{L^2}\leq
    \norm{ f(x)^{\frac{2}{3}}u}_{L^2}\leq  C\set{\norm{ Pu}_{L^2}+
    \norm{u}_{L^2}}.
  \end{eqnarray}
  The above two inequalities yield that for any $u\in C_0^\infty\inner{\mathbb{R}^{2n}},$
  \begin{eqnarray}\label{7281}
    \norm{\Lambda_y  f(x)^{s\over2}u}_{L^2}^2
    \leq C \set{\abs{\comi{
    P f(x)^{s\over2}u,~ f(x)^{s\over2}u}_{L^2}}
    +\norm{Pu}_{L^2}^2+\norm{u}_{L^2}^2}.
  \end{eqnarray}

  \textbf{a)} ~Firstly let us consider the case when
   $s\leq \frac{2}{3}.$ In such a case, we have
   \bigskip
  \begin{eqnarray*}
    \abs{\comi{ P f(x)^{s\over2}u,~ f(x)^{s\over2}u}_{L^2}}&\leq&
    \abs{\comi{ P u,~ f(x)^{s}u}_{L^2}}+\abs{\comi{ \Big[P,~ f(x)^{s\over2}\Big]u,~
    f(x)^{s\over2}u}_{L^2}}\\
    &\leq&
    \norm{Pu}_{L^2}^2+\norm{f(x)^{2\over3}u}_{L^2}^2+\abs{\comi{ \Big[P,~ f(x)^{s\over2}\Big]u,~
    f(x)^{s\over2}u}_{L^2}}\\
    &\leq&
    C\norm{Pu}_{L^2}^2+C\norm{u}_{L^2}^2+\abs{\comi{ \Big[P,~ f(x)^{s\over2}\Big]u,~
    f(x)^{s\over2}u}_{L^2}},
  \end{eqnarray*}
  the last inequality using \reff{req}. This along with \reff{7281}
  gives
   \begin{eqnarray*}
    \norm{\Lambda_y  f(x)^{s\over2}u}_{L^2}^2
    \leq C\norm{Pu}_{L^2}^2+C\norm{u}_{L^2}^2+\abs{\comi{ \Big[P,~ f(x)^{s\over2}\Big]u,~
    f(x)^{s\over2}u}_{L^2}}.
  \end{eqnarray*}
  On the other hand using \reff{091005+1}  with $s\leq {2\over3}$
  implies, for any $\eps>0,$
  \begin{eqnarray*}
    \abs{\comi{ \Big[P,~ f(x)^{s\over2}\Big]u,~
    f(x)^{s\over2}u}_{L^2}}\leq C\norm{\Lambda_yf(x)^{s\over2}u}_{L^2}\,\norm{
    u}_{L^2}\leq \eps\norm{\Lambda_yf(x)^{s\over2}u}_{L^2}^2+C_\eps\norm{
    u}_{L^2}^2.
  \end{eqnarray*}
  Combining the above two inequalities and taking $\eps\leq {1\over2}$, we get
  \begin{eqnarray*}
    \norm{\Lambda_y  f(x)^{s\over2}u}_{L^2}^2
    \leq C\norm{Pu}_{L^2}^2+C\norm{u}_{L^2}^2.
  \end{eqnarray*}
   Since $\sum\limits_{\mu\geq1}\norm{R_\mu
    u}_{L^2}^2\leq C
    \norm{\Lambda_y  f(x)^{s\over2}u}_{L^2}^2,$ then the above estimate gives the validity of
    \reff{CMU}
    for $s\leq {2\over3}$.

\textbf{b)}~Next we shall prove \reff{CMU} for
${2\over3}<s<{4\over3}$. If $\frac{10}{9}<s<\frac{4}{3},$ then it
follows from \reff{7281} and \reff{0909261} that
\begin{eqnarray}\label{++ad1}
   \forall  C_0^\infty\inner{\mathbb{R}^{2n}},\quad  \norm{\Lambda_y  f(x)^{s\over2}u}_{L^2}^2
    \leq C\norm{
    P f(x)^{s\over2}u}_{L^2}^2+\norm{ Pu}_{L^2}
    +C\norm{u}_{L^2}^2.
\end{eqnarray}
This gives the validity of \reff{CMU} for $s\in
]\frac{10}{9},~\frac{4}{3}[.$

  Now we focus on the case when $\frac{2}{3}<s\leq \frac{10}{9}.$
  Observe that
  \begin{eqnarray*}
    &\abs{\comi{
    P f(x)^{s\over2}u,~ f(x)^{s\over2}u}_{L^2}}=
    \abs{\comi{
    P f(x)^{s\over2}u,~ f(x)^{{2\over3}+\inner{\frac{s}{2}-{2\over3}}}u}_{L^2}} \\
    &\leq
    \abs{\comi{
    P f(x)^{s-{2\over3}}u,~ f(x)^{{2\over3}}u}_{L^2}}
    +\abs{\comi{\com{
    P,~ f(x)^{\frac{s}{2}-{2\over3}}} f(x)^{{s\over2}}u,
    ~ f(x)^{{2\over3}}u}_{L^2}}.
  \end{eqnarray*}
  Moreover since $f(x)$ satisfies \reff{091005+1}, then
  \begin{eqnarray*}
    \abs{\com{
    P,~ f(x)^{\frac{s}{2}-{2\over3}}} f(x)^{{s\over2}}u}
    &\leq C\abs{y} f(x)^{2s-{5\over3}}\abs{u},
  \end{eqnarray*}
  and thus
  \[\abs{\comi{\com{
    P,~ f(x)^{\frac{s}{2}-{2\over3}}} f(x)^{{s\over2}}u,
    ~ f(x)^{{2\over3}}u}_{L^2}}\leq \eps \norm{\Lambda_y f(x)^{2s-{5\over3}}u}_{L^2}^2
    +C_\eps
    \norm{f(x)^{{2\over3}}u}_{L^2}^2
    \]
  Combination of the above three inequalities gives
  \begin{eqnarray*}
    \abs{\comi{
    P f(x)^{s\over2}u,~ f(x)^{s\over2}u}_{L^2}}
    \leq\eps \norm{\Lambda_y f(x)^{2s-{5\over3}}u}_{L^2}^2+
    C_\eps\set{\norm{
    P f(x)^{s-{2\over3}}u}_{L^2}^2 +\norm{
    f(x)^{{2\over3}}u}_{L^2}^2}.
  \end{eqnarray*}
  Moreover since $2s-{5\over3}\leq \frac{s}{2}$ for $s\leq\frac{10}{9},$ then
   \begin{eqnarray*}
    \norm{\Lambda_y f(x)^{2s-{5\over3}}u}_{L^2}^2
    \leq
    \norm{\Lambda_y f(x)^{{s\over2}}u}_{L^2}^2,
   \end{eqnarray*}
   and hence by \reff{req} we obtain
    \begin{eqnarray*}
    \abs{\comi{
    P f(x)^{s\over2}u,~ f(x)^{s\over2}u}_{L^2}}
    \leq \eps \norm{\Lambda_y f(x)^{{s\over2}}u}_{L^2}^2+
    C_{\eps}\set{\norm{
    P f(x)^{s-{2\over3}}u}_{L^2}^2+\norm{
    Pu}_{L^2}^2+\norm{u}_{L^2}^2}.
  \end{eqnarray*}
   Inserting the above inequality into \reff{7281} and then taking $\eps$ small enough,
   we get the desired estimate \reff{CMU} for
  ${2\over3}<s<{10\over9}.$ Thus the proof of Lemma \ref{patch} is completed.
  \qed

\medskip

  Now we are ready to prove the main result of this section.

  \noindent\emph{Proof of Proposition \ref{prpfin}.}
   Now we want to show that
  \begin{eqnarray}\label{+09865++}
    \norm{\inner{1-\triangle_x}^{{\delta\over2}}u}_{L^2}^2\leq
    C\set{\norm{P u}_{L^2}^2
    +\norm{u}_{L^2}^2}.
  \end{eqnarray}
  Recall $\delta$  equals to $\frac{2}{3}$ if $s\leq \frac{2}{3}$,
  $\frac{4}{3}-s$ if ${2\over3}<s\leq \frac{10}{9},$ and
  $\frac{2}{3}-{s\over2}$ if ${{10}\over9}< s<\frac{4}{3}.$
  Using the estimates \reff{expression} and \reff{CMU} gives
  that
  \begin{eqnarray}\label{091005+2}
    \forall~u\in C_0^\infty\big(\mathbb{R}^{2n}\big),
    \quad \sum_{\mu\geq1}
    \norm{P_{x_\mu} \varphi_\mu\, u}_{L^2}^2
    \leq C
    \set{\norm{P f(x)^{\tilde s}
    u}_{L^2}^2+\norm{P
    u}_{L^2}^2+\norm{u}_{L^2}^2},
  \end{eqnarray}
  where $\tilde s=\frac{2}{3}-\delta.$ We can verify that
  \begin{eqnarray}\label{091005+3}
    -\tilde s+s-1\leq 0.
  \end{eqnarray}

  Firstly let us consider the case of $s\leq
  \frac{2}{3}.$  Then $\tilde s=0$ and \reff{091005+2} becomes
    \[\forall~u\in C_0^\infty\big(\mathbb{R}^{2n}\big),
    \quad \sum_{\mu\geq1}
    \norm{P_{x_\mu} \varphi_\mu  u}_{L^2}^2
    \leq C
    \set{\norm{P
    u}_{L^2}^2+\norm{u}_{L^2}^2}.
  \]
  On the other hand, using \reff{tra} with $a=\frac{1}{3}$ and then \reff{+ES}, we have
   \begin{eqnarray*}
   \norm{\inner{1-\triangle_x}^{{1\over3}}u}_{L^2}^2
   &\leq& C
   \sum_{\mu\geq1}\norm{\inner{1-\triangle_x}^{{1\over3}}\varphi_\mu\,
   u}_{L^2}^2+C\norm{Pu}_{L^2}^2+C\norm{u}_{L^2}^2\\&\leq& C \sum_{\mu\geq1}
   \norm{P_{x_\mu} \varphi_\mu\, u}_{L^2}^2
   +C\sum_{\mu\geq1}\norm{\varphi_\mu\,
   u}_{L^2}^2+C\norm{Pu}_{L^2}^2+C\norm{u}_{L^2}^2.
  \end{eqnarray*}
  As a result, we get from these inequalities
  \begin{eqnarray*}
    \forall~u\in C_0^\infty\big(\mathbb{R}^{2n}\big)
    \quad \norm{
    \inner{1-\triangle_x}^{1\over3}u}_{L^2}^2
    \leq C
    \set{\norm{P
    u}_{L^2}^2+\norm{u}_{L^2}^2}.
  \end{eqnarray*}
  This gives the validity of \reff{+09865++} for $s\leq {2\over3}.$

  Now we consider the case when
  $\frac{2}{3}<s<\frac{4}{3}.$
  Note that $\delta=\frac{2}{3}-\tilde s.$  Then we use \reff{tra} with $a=\frac{\delta}{2}$ to
  get
  \begin{equation*}
  \begin{split}
    &\norm{\inner{1-\triangle_x}^{{\delta\over2}}u}_{L^2}^2
    \leq  C\sum_{\mu\geq
   1}\norm{\inner{1-\triangle_x}^{\frac{\delta}{2}}\varphi_\mu
    u}_{L^2}^2+C\norm{Pu}_{L^2}^2+C\norm{u}_{L^2}^2\\
    &\leq C\sum_{\mu\geq
    1}\norm{\Big(1+\frac{1}{2}\abs{\partial_xV(x_\mu)}^2-\triangle_x\Big)^{\delta\over2}\varphi_\mu
    u}_{L^2}^2+C\norm{Pu}_{L^2}^2+C\norm{u}_{L^2}^2\\
    &= C\sum_{\mu\geq
    1}\norm{\Big(1+\frac{1}{2}\abs{\partial_xV(x_\mu)}^2-\triangle_x\Big)^{1\over3}
    \Big(1+\frac{1}{2}\abs{\partial_xV(x_\mu)}^2-\triangle_x\Big)^{-\frac{\tilde
    s}{2}}\varphi_\mu
    u}_{L^2}^2\\
    &\quad+C\norm{Pu}_{L^2}^2+C\norm{u}_{L^2}^2\\
    &\leq C\sum_{\mu\geq
    1}\norm{\Big(1+\frac{1}{2}\abs{\partial_xV(x_\mu)}^2-\triangle_x\Big)^{1\over3}
    f(x_\mu)^{-\tilde s}\varphi_\mu
    u}_{L^2}^2+C\norm{Pu}_{L^2}^2+C\norm{u}_{L^2}^2.
  \end{split}
  \end{equation*}
  Consequently, using \reff{ES} yields
  \begin{eqnarray*}
    \norm{\inner{1-\triangle_x}^{{\delta\over2}}u}_{L^2}^2
    \leq C\sum_{\mu\geq 1}\norm{P_{x_\mu} ~f(x_\mu)^{-\tilde s}\varphi_\mu\, u}_{L^2}^2
    +C\norm{Pu}_{L^2}^2+C\norm{u}_{L^2}^2.
  \end{eqnarray*}
  Thus \reff{+09865++} will follow if we can show that
  \begin{eqnarray}\label{09865}
    \sum_{\mu\geq 1}\norm{P_{x_\mu} ~f(x_\mu)^{-\tilde s} \varphi_\mu\,u}_{L^2}^2
    \leq C\set{\norm{P u}_{L^2}^2
    +\norm{u}_{L^2}^2}.
  \end{eqnarray}
  To prove \reff{09865}, we write
  \[
    f(x_\mu)^{-\tilde s}\varphi_\mu
    =\inner{f(x)^{\tilde s}f(x_\mu)^{-\tilde s}}
    \varphi_\mu  \,f(x)^{-\tilde s}.
  \]
  Then
  \begin{eqnarray*}
    \sum_{\mu\geq 1}\norm{P_{x_\mu} ~f(x_\mu)^{-\tilde s}\varphi_\mu\, u}_{L^2}^2
    \leq (I)+(II)
  \end{eqnarray*}
  with $(I),(II)$ given by
  \[
    (I)=2 \sum_{\mu\geq 1}\norm{\inner{ f(x)^{\tilde s}f(x_\mu)^{-\tilde
    s}}
    P_{x_\mu} ~\varphi_\mu   f(x)^{-\tilde s}
    u}_{L^2}^2
  \]
  and
  \[
    (II)=2\sum_{\mu\geq 1}\norm{\com{P_{x_\mu},~ f(x)^{\tilde s}f(x_\mu)^{-\tilde s}
    }\varphi_\mu   f(x)^{-\tilde s}
    u}_{L^2}^2.
  \]
  By \reff{7310}, we see
  \begin{eqnarray*}
    (I)\leq C\sum_{\mu\geq 1}\norm{
    P_{x_\mu} ~\varphi_\mu   f(x)^{-\tilde s}
    u}_{L^2}^2.
  \end{eqnarray*}
  This along with \reff{091005+2} gives
  \begin{eqnarray}\label{+091005}
    (I)\leq C
    \set{\norm{P
    u}_{L^2}^2+\norm{P
     f(x)^{-\tilde s}u}_{L^2}^2+\norm{ f(x)^{-\tilde
    s}u}_{L^2}^2}.
  \end{eqnarray}
  By use of \reff{091005+1} and \reff{091005+3}, we
  have
  \begin{eqnarray*}
    \norm{\big[P,~ f(x)^{-\tilde s}\,\big]u}_{L^2}^2
    \leq  C \norm{ f(x)^{-\tilde s+s-1} \abs{y} u}_{L^2}^2
    \leq  C \norm{\abs{y} u}_{L^2}^2,
  \end{eqnarray*}
  and hence
  \begin{eqnarray*}
    \norm{P f(x)^{-\tilde s}u
    }_{L^2}^2
    \leq
    2\norm{Pu}_{L^2}^2
    +2\norm{\big[P,~ f(x)^{-\tilde s}\,\big]u}_{L^2}^2\leq C\set{\norm{P u}_{L^2}^2+\norm{
    u}_{L^2}^2}.
  \end{eqnarray*}
  This along with \reff{+091005} gives
  \begin{eqnarray*}
    I\leq C\set{\norm{P u}_{L^2}^2+\norm{ u}_{L^2}^2}.
  \end{eqnarray*}
  Now it remains to treat the term $(II).$ The equality
  \[
    \com{P_{x_\mu},~f(x)^{\tilde s}f(x_\mu)^{-\tilde
    s}}= \inner{y\cdot\partial_x \inner{f(x)^{\tilde s}}}
    f(x_\mu)^{-\tilde s}
  \]
  gives
  \begin{eqnarray}\label{Fin}
    (II)=2\sum_{\mu\geq 1}\norm{\inner{y\cdot\partial_x \inner{f(x)^{\tilde
    s}}}
    f(x_\mu)^{-\tilde s}
    f(x)^{-\tilde s}\varphi_\mu
    u}_{L^2}^2.
  \end{eqnarray}
 By \reff{7310}, \reff{091005+1}  and \reff{091005+3},  we have
  \begin{eqnarray*}
    \abs{\partial_x \inner{f(x)^{\tilde s}}}f(x_\mu)^{-\tilde s}
     f(x)^{-\tilde s}\varphi_\mu\leq
    C f(x)^{s-1-\tilde s}\leq C.
  \end{eqnarray*}
  So
  \begin{eqnarray*}
    \sum_{\mu\geq 1}\norm{\inner{y\cdot\partial_x \inner{f(x)^{\tilde
    s}}}
    f(x_\mu)^{-\tilde s}
    f(x)^{-\tilde s}\varphi_\mu\,
    u}_{L^2}^2\leq C
    \sum_{\mu\geq 1}\norm{\varphi_\mu\abs{y}
    u}_{L^2}^2
    \leq C\norm{\Lambda_yu}_{L^2}^2.
  \end{eqnarray*}
  This along with \reff{Fin} gives
  \begin{eqnarray*}
    (II)\leq C\norm{\Lambda_yu}_{L^2}^2\leq C\set{\norm{Pu}_{L^2}^2
    +\norm{u}_{L^2}^2}.
  \end{eqnarray*}
  Combining the estimate on the term $(I),$ we get the required
  inequality \reff{09865}.  The proof of Proposition \ref{prpfin}
  is thus completed.
  \qed

\section{Proof of Corollary \ref{thm+2}}\label{sec5}

The method is quite similar as that in \cite{HelfferNier05}, and the
main difference is that we have to use the functional calculus for
self-adjoint operators instead of the pseudo-differential calculus
used in \cite{HelfferNier05}, since in our case the potential $V$
only belongs to $C^2\big(\mathbb{R}^{2n}\big)$.  Firstly let us
mention some well-known facts on the functional calculus for
positive self-adjoint operators (see for instance Chapter XI of
\cite{Yosida95}). Consider the  Schr\"{o}dinger operator
$1+\abs{\partial_x V(x)}^2-\triangle_x$ which is defined on
$C_0^\infty(\mathbb{R}^n).$ Since
\[
1+\abs{\partial_x V(x)}^2-\triangle_x\geq 1
\]
holds in the sense of operators, then it is well-known that
$1+\abs{\partial_x V(x)}^2-\triangle_x$ admits a unique self-adjoint
extension on $L^2(\mathbb{R}^{n})$, still denoted by
$1+\abs{\partial_x V(x)}^2-\triangle_x.$ As a result, using the
notation
\[
A=1+\abs{\partial_x V(x)}^2-\triangle_x,
\]
the self-adjoint operator $A$ admits a spectral representation
\[
A=\int_0^{+\infty} \lambda dE_\lambda
\]
with domain
\[
  D(A)=\set{u\in L^2(\mathbb{R}^n);\quad \int_0^{+\infty} \lambda^2
  d\,\|E_\lambda u\|_{L^2}^2<+\infty}.
\]
Here $\set{E_\lambda}_{\lambda\geq0}$ is called a spectral
resolution of $1+\abs{\partial_x V(x)}^2-\triangle_x$. By the
spectral representation, we can define the fractional power  of the
operator $A$ as follows: for each $\theta\geq 0,$
\[
A^\theta=\inner{1+\abs{\partial_x V(x)}^2-\triangle_x}^\theta
=\int_0^{+\infty} \lambda^\theta dE_\lambda
\]
with domain
\[\quad D(A^\theta)=
\set{u\in L^2(\mathbb{R}^n);\quad \int_0^{+\infty} \lambda^{2\theta}
  d\,\|E_\lambda u\|_{L^2}^2<+\infty}.
\]
Note that $\inner{1+\abs{\partial_x V(x)}^2-\triangle_x}^0=I$.
Moreover, since $A=1+\abs{\partial_x V(x)}^2-\triangle_x\geq 1,$
then $E_\lambda=0$ for $0\leq \lambda<1.$ This allows us to define
the negative fraction power by
\[
A^{-\theta} =\int_0^{+\infty} \lambda^{-\theta}
dE_\lambda,\quad\theta>0,
\]
with domain
\[\quad D(A^{-\theta})=
\set{u\in L^2(\mathbb{R}^n);\quad \int_0^{+\infty}
\lambda^{-2\theta}
  d\,\|E_\lambda u\|_{L^2}^2<+\infty}.
\]
 Now we list some classical results to be used frequently on the
fractional power of the operator $A$. For each $\theta\geq 0,$ the
operators $A^{\pm \theta}$ are self-adjoint on
$L^2\big(\mathbb{R}^{n}\big),$  and satisfy the following relation
\[
    \inner{1+\abs{\partial_x V(x)}^2-\triangle_x}^{-\theta}=
    \inner{\inner{1+\abs{\partial_x
    V(x)}^2-\triangle_x}^{\theta}}^{-1},
\]
that is, $A^{-\theta}$ is the inverse operator of $A^{\theta}$. If
$\theta\in[0,1]$ then $u\in D(A)$ if and only if $u\in D(A^\theta)$
and also $A^\theta\in D(A^{1-\theta}).$ For such $u$ we have
\[
Au=A^{1-\theta}A^\theta u.
\]
Moreover the negative power $A^{-\theta},\theta>0,$ is bounded on
$L^2(\mathbb{R}^n),$ and for any $\theta_1,\theta_2\geq 0,$ and any
$u\in L^2\big(\mathbb{R}^{n}\big),$ we have
\[
      A^{-(\theta_1+\theta_2)}u=
      A^{-\theta_1} A^{-\theta_2}u.
\]
Since the following inequalities
\[
1+\abs{\partial_x V(x)}^2-\triangle_x\geq 1+\abs{\partial_x
V(x)}^2~~\,\, {\rm and}~~ ~~ 1+\abs{\partial_x V(x)}^2-\triangle_x
\geq 1-\triangle_x
\]
hold in the sense of operators, then the monotonicity of operator
functional implies for each $\theta\geq0,$
\[
\inner{1+\abs{\partial_x V(x)}^{2}-\triangle_x}^{-2\theta}\leq
\inner{1+\abs{\partial_x V(x)}^2}^{-2\theta}
\]
and
\[
\inner{1+\abs{\partial_x V(x)}^2-\triangle_x}^{-2\theta} \leq
\inner{1-\triangle_x}^{-2\theta}.
\]
As a result, for each $\theta\geq0$ the following estimates hold:
for any $u\in L^2\big(\mathbb{R}^{2n}\big),$
\begin{eqnarray}\label{+es++1}
  \norm{\inner{1+\abs{\partial_x V(x)}^{2}-\triangle_x}^{-\theta}u}_{L^2}\leq
  \norm{\inner{1+\abs{\partial_x V(x)}^2}^{-\theta}u}_{L^2}
\end{eqnarray}
and
\begin{eqnarray}\label{+es++2}
  \norm{\inner{1+\abs{\partial_x V(x)}^{2}-\triangle_x}^{-\theta}u}_{L^2}\leq
  \norm{\inner{1-\triangle_x}^{-\theta}u}_{L^2}.
\end{eqnarray}

\bigskip
The rest of the paper is devoted to showing  Corollary \ref{thm+2}.

\noindent\emph{Proof  of Corollary \ref{thm+2}.} Assume $V\in
C^2\big(\mathbb{R}^{2n}\big)$ satisfies \reff{+HyP}. Then the global
hypoelliptic estimates \reff{+A1} and \reff{+++A1} hold. Let
$\delta$ be the number given in \reff{+++A1}. For any $u\in
C_0^\infty(\mathbb{R}^{2n}),$ we can verify that
\begin{eqnarray*}
  &&\norm{\inner{1+\abs{\partial_xV(x)}^2-\triangle_x}^{\delta\over2}u}_{L^2}=
  \norm{\inner{1+\abs{\partial_xV(x)}^2-\triangle_x}^{{\delta\over2}-1}
  \inner{1+\abs{\partial_xV(x)}^2-\triangle_x}u}_{L^2}\\
  &&\leq\norm{\inner{1+\abs{\partial_xV(x)}^2-\triangle_x}^{{\delta\over2}-1}
  \abs{\partial_xV(x)}^2u}_{L^2}+
  \norm{\inner{1+\abs{\partial_xV(x)}^2-\triangle_x}^{{\delta\over2}-1}
  \inner{1-\triangle_x}u}_{L^2}\\
  &&\leq\norm{
  \abs{\partial_xV(x)}^{\delta}u}_{L^2}+
  \norm{\inner{1-\triangle_x}^{{\delta\over2}}u}_{L^2}+\norm{u}_{L^2},
\end{eqnarray*}
the last inequality using from \reff{+es++1} and \reff{+es++2} with
$-\theta={\delta\over2}-1<0$. As a result, applying \reff{+A1} and
\reff{+++A1} yields
\begin{eqnarray*}
  \forall~u\in C_0^\infty\big(\mathbb{R}^{2n}\big),\quad \norm{\inner{1+\abs{\partial_xV(x)}^2-\triangle_x}^{\delta\over2}u}_{L^2}\leq
  C\set{\norm{Pu}_{L^2}+\norm{u}_{L^2}}.
\end{eqnarray*}
This implies the operator
\[
 \inner{1+\abs{\partial_xV(x)}^2-\triangle_x}^{\delta\over2}\circ\inner{1+P}^{-1}:~~L^2\big(\mathbb{R}^{2n}\big)\longrightarrow
 L^2\big(\mathbb{R}^{2n}\big)
\]
is bounded. Since
\[
\inner{1+P}^{-1}=\inner{1+\abs{\partial_xV(x)}^2-\triangle_x}^{-{\delta\over2}}\circ\inner{
\inner{1+\abs{\partial_xV(x)}^2-\triangle_x}^{\delta\over2}\circ\inner{1+P}^{-1}},
\]
then the compactness of the resolvent $\inner{1+P}^{-1}$  will
follow if the operator
\[
\inner{1+\abs{\partial_xV(x)}^2-\triangle_x}^{-{\delta\over2}}:~~L^2\big(\mathbb{R}^{2n}\big)\longrightarrow
 L^2\big(\mathbb{R}^{2n}\big)
\]
is compact. Recall that a operator $B$ acting on
$L^2\big(\mathbb{R}^{2n}\big)$ is compact if and only if
\[
\forall~u_{n}\in L^2\big(\mathbb{R}^{2n}\big),~ u_n\stackrel{{\rm
weakly}}{\longrightarrow} 0\Longrightarrow
\lim_{n\rightarrow+\infty}\norm{Bu_n}_{L^2}=0.
\]
Then to get the compactness of
$\inner{1+\abs{\partial_xV(x)}^2-\triangle_x}^{-{\delta\over2}}$, we
have to show that for any sequence $\set{u_n}_{n\geq1}$ converging
to $0$ weakly in $L^2(\mathbb{R}^{2n})$,
\begin{eqnarray}\label{0910012}
  \lim_{n\rightarrow+\infty}\norm{\inner{1+\abs{\partial_xV(x)}^2-\triangle_x}^{-{\delta\over2}}u_n}_{L^2}=0.
\end{eqnarray}
This can be derived from the compactness of resolvent of the Witten
Laplacian $\triangle_{V/2}^{(0)}.$ Indeed, since
$\triangle_{V/2}^{(0)}$ has a compact resolvent then the operator
$\inner{1+\abs{\partial_xV(x)}^2-\triangle_x}^{-1}$ is  also
compact. As a result,
\begin{eqnarray*}
\lim_{n\rightarrow+\infty}\norm{\inner{1+\abs{\partial_xV(x)}^2-\triangle_x}^{-1}u_n}_{L^2}=0.
\end{eqnarray*}
Moreover since $\set{u_n}_{n\geq1}$ is a weakly convergent sequence
then it is bounded in $L^2\big(\mathbb{R}^{2n}\big)$. Hence
\begin{eqnarray}\label{0910011}
\lim_{n\rightarrow+\infty}\norm{\inner{1+\abs{\partial_xV(x)}^2-\triangle_x}^{-1}u_n}_{L^2}\norm{u_n}_{L^2}=0.
\end{eqnarray}
As a result, using the relation
\begin{eqnarray*}
  \norm{\inner{1+\abs{\partial_xV(x)}^2-\triangle_x}^{-{1\over2}}u_n}_{L^2}^2
  =\comi{\inner{1+\abs{\partial_xV(x)}^2-\triangle_x}^{-1}u_n,
  ~u_n}_{L^2}
\end{eqnarray*}
gives
\begin{eqnarray*}
  \lim_{n\rightarrow+\infty}\norm{\inner{1+\abs{\partial_xV(x)}^2-\triangle_x}^{-{1\over2}}u_n}_{L^2}=0,
\end{eqnarray*}
and hence
\begin{eqnarray*}
  \lim_{n\rightarrow+\infty}\norm{\inner{1+\abs{\partial_xV(x)}^2-\triangle_x}^{-{1\over2}}u_n}_{L^2}\norm{u_n}_{L^2}=0.
\end{eqnarray*}
Repeating the above arguments,  we can get
\begin{eqnarray*}
  \lim_{n\rightarrow+\infty}\norm{\inner{1+\abs{\partial_xV(x)}^2-\triangle_x}^{-{\delta\over2}}u_n}_{L^2}=0.
\end{eqnarray*}
Then the proof of Corollary \ref{thm+2}  is completed. \qed

\section{Appendix}
Here we present another proof of Proposition \ref{PA}. Let's restate
it as
\begin{proposition}\label{++PA+} We associate with each
fixed $x_0\in\mathbb{R}^{n}$ the operator
\[
  P_{x_0}=y \cdot \partial_x -\partial_x V(x_0)\cdot
  \partial_y-\triangle_y+\frac{\abs y^2}{4}-\frac{n}{2}=y \cdot \partial_x -\partial_x V(x_0)\cdot
  \partial_y+\sum_{j=1}^n
  L_j^*L_j.
\]
Then there is a constant $C$ independent of $x_0,$  such that
 for any $u\in  C_0^\infty\big(\mathbb{R}^{2n}\big),$ one has
  \begin{eqnarray}\label{+ES2010}
   \abs{\partial_x
   V(x_0)}^{4\over3}\norm{ u}_{L^2}^2+\norm{\inner{1-\triangle_x}^{1\over3} u}_{L^2}^2
   \leq
   C\set{
   \norm{P_{x_0} u}_{L^2}^2+\norm{u}_{L^2}^2},
\end{eqnarray}
or equivalently,
\begin{eqnarray}\label{ES2010}
   \norm{\tilde \Lambda_{x_0}^{2\over3}u}_{L^2}^2 \leq
   C\set{
   \norm{P_{x_0} u}_{L^2}^2+\norm{u}_{L^2}^2},
\end{eqnarray}
where $\tilde \Lambda_{x_0}=\inner{1+\frac{1}{2}\abs{\partial_x
   V(x_0)}^2-\triangle_x}^{1\over2}.$

\end{proposition}

\noindent\emph{Proof.} We will prove \reff{ES2010} in this proof. To
simplify the notation the capital letter $C$ will be used to denote
different suitable constants independent of $x_0.$ Denoting
$Q_{x_0}=y \cdot D_x-\partial_x V(x_0)\cdot
 D_y,$ we can write
$P_{x_0}$ as
\begin{eqnarray*}
  P_{x_0}=iQ_{x_0}+\sum_{j=1}^nL_j^*L_j,
\end{eqnarray*}
from which we deduce
\begin{eqnarray}\label{Non}
  \forall~u\in C_0^\infty\big(\mathbb{R}^{2n}\big),\quad \sum_{j=1}^n\norm{L_j u}_{L^2}^2\leq {\rm Re}\comi{P_{x_0}u,~u}_{L^2}  .
\end{eqnarray}
Then
\begin{eqnarray}\label{8501}
  \norm{\Lambda_yu}_{L^2}^2\leq C\set{\sum_{j=1}^n\norm{L_j
  u}_{L^2}^2+\norm{u}_{L^2}^2}
  \leq C\set{{\rm Re}\comi{P_{x_0}u,~u}_{L^2}+\norm{u}_{L^2}^2}.
\end{eqnarray}
We will  proceed to prove \reff{ES2010} in the following three
steps.

 {\bf\emph{Step A}.}~ We claim, for any $\eps>0,$ there is a constant $C_{\eps},$ depending only on $\eps,$
 such that
 \begin{eqnarray}\label{LJ}
  \sum_{j=1}^n\norm{L_j^*L_j u}_{L^2}^2\leq\eps\norm{\tilde \Lambda_{x_0}^{2/3}u}_{L^2}^2+C_\eps\set{\norm{P_{x_0}u}_{L^2}^2
   +\norm{u}_{L^2}^2}.
\end{eqnarray}
To confirm this, we apply  \reff{Non} to get, for any $u\in
 C_0^\infty\big(\mathbb{R}^{2n}\big),$
\begin{eqnarray*}
  \norm{L_jL_j^*u}_{L^2}^2
  &\leq& {\rm Re}\comi{P_{x_0}L_j^*u,~L_j^*
  u}_{L^2}\\
  &=&{\rm Re}\comi{[P_{x_0},~L_j^*]u,~L_j^*
  u}_{L^2}+ {\rm Re}\comi{P_{x_0}u,~L_jL_j^*
  u}_{L^2}\\
  &\leq&{\rm Re}\comi{[P_{x_0},~L_j^*]u,~L_j^*
  u}_{L^2}+{1\over2}
  \norm{L_jL_j^*u}_{L^2}^2+2\norm{P_{x_0}u}_{L^2}^2.
\end{eqnarray*}
Hence
\begin{eqnarray}\label{090922}
\norm{L_jL_j^*u}_{L^2}^2
  \leq 2{\rm Re}\comi{[P_{x_0},~L_j^*]u,~L_j^*
  u}_{L^2}+4\norm{P_{x_0}u}_{L^2}^2.
\end{eqnarray}
To estimate the first term on the right side of the above
inequality, we make use of the following commutation relations
satisfied by $iQ_{x_0}, L_j, L_k^*, j,k=1,2,\cdots,n,$
\[
  [iQ_{x_0},L_j^*]=-\frac{1}{2}\partial_{x_j}
   V(x_0)+\partial_{x_j},\quad [L_j,~L_k]=[L_j^*,~L_k^*]=0,\quad [L_j,~L_k^*]=\delta_{jk};
\]
this gives, for any $\tilde\eps>0,$
\begin{eqnarray*}
   {\rm Re}\comi{[P_{x_0},~L_j^*]u,~L_j^*u}_{L^2}&=&
   \comi{L_j^*u,~L_j^*u}_{L^2}
   +\Big<\Big(-\frac{1}{2}\partial_{x_j}
   V(x_0)+\partial_{x_j}\Big)u,~L_j^*u\Big>_{L^2}\\
   &\leq& \tilde\eps\norm{\tilde\Lambda_{x_0}^{2\over3}~u}_{L^2}^2
   +C_{\tilde\eps}\set{\norm{L_j\tilde\Lambda_{x_0}^{1\over3}u}_{L^2}^2+
   \norm{L_ju}_{L^2}^2+\norm{u}_{L^2}^2}.
\end{eqnarray*}
Moreover we use \reff{Non} again to obtain
\begin{eqnarray*}
   \norm{L_j\tilde\Lambda_{x_0}^{1\over3}u}_{L^2}^2
   \leq{\rm Re}\comi{P_{x_0}\tilde\Lambda_{x_0}^{1\over3}u,~\tilde\Lambda_{x_0}^{1\over3}
  u}_{L^2}={\rm Re}\comi{P_{x_0}u,~\tilde\Lambda_{x_0}^{2\over3}
  u}_{L^2}.
\end{eqnarray*}
Combining these  inequalities, we conclude
\begin{eqnarray*}
   {\rm Re}\comi{[P_{x_0},~L_j^*]u,~L_j^*u}_{L^2}
   \leq \eps\norm{\tilde\Lambda_{x_0}^{2\over3}~u}_{L^2}^2
   +C_{\eps}\set{\norm{P_{x_0}u}_{L^2}^2+
   \norm{u}_{L^2}^2},
\end{eqnarray*}
and thus by \reff{090922}
\begin{eqnarray*}
   \norm{L_jL_j^*u}_{L^2}^2
   \leq \eps\norm{\tilde\Lambda_{x_0}^{2\over3}~u}_{L^2}^2
   +C_{\eps}\set{\norm{P_{x_0}u}_{L^2}^2+
   \norm{u}_{L^2}^2}.
\end{eqnarray*}
Observe
\[
   \norm{L_j^*L_j u}_{L^2}^2\leq 2\norm{L_jL_j^*
   u}_{L^2}^2+2\norm{u}_{L^2}^2.
\]
Then the desired estimate \reff{LJ} follows.

{\bf \emph{Step B}.}~  For $x_0\in\mathbb{R}^n$  consider the
operator
\[
   \ell_{x_0}= \tilde\Lambda_{x_0}^{-{2\over3}}\circ\set{\partial_x V(x_0)\cdot y+2D_x\cdot D_y}.
\]
It's a straightforward verification to see that
\begin{eqnarray}\label{Qell}
  {\rm Re}\comi{iQ_{x_0}u,~\ell_{x_0}u}_{L^2}=
  \frac{i}{2}\comi{\com{\ell_{x_0},~Q_{x_0}}u,~u}_{L^2}
  =\norm{\tilde\Lambda_{x_0}^{2/3}u}_{L^2}^2-\norm{\tilde\Lambda_{x_0}^{-1/3}~u}_{L^2}^2.
\end{eqnarray}
Next we will  prove that, for any $\eps>0,$
  \begin{eqnarray}\label{L2}
    \norm{\ell_{x_0} u}_{L^2}^2\leq \eps\norm{\tilde\Lambda_{x_0}^{2/3}u}_{L^2}^2+ C_\eps
    \set{\norm{P_{x_0}u}_{L^2}^2+\norm{u}_{L^2}^2}.
  \end{eqnarray}
To confirm this,  observe
\[
  \norm{\ell_{x_0} u}_{L^2}^2\leq 2\norm{\tilde\Lambda_{x_0}^{-2/3}~\partial_{x}V(x_0)\cdot y~ u}_{L^2}^2+
   4\norm{\Lambda_{x_0}^{-2/3}D_x\cdot D_y u}_{L^2}^2
   \leq C\norm{\tilde\Lambda_{x_0}^{1/3}  ~\Lambda_y u}_{L^2}^2.
\]
This along with \reff{8501} gives that
\begin{eqnarray*}
  \norm{\ell_{x_0} u}_{L^2}^2&\leq& C\norm{\Lambda_{x_0}^{1/3}  ~\Lambda_y u}_{L^2}^2\leq C
  \set{{\rm Re}\comi{P_{x_0}\tilde\Lambda_{x_0}^{1/3} u,~\tilde\Lambda_{x_0}^{1/3}u}_{L^2}+\norm{\tilde\Lambda_{x_0}^{1/3}u}_{L^2}^2}
  \\&\leq& C
  \set{\abs{\comi{P_{x_0}
  u,~\tilde\Lambda_{x_0}^{2/3}u}_{L^2}}+\norm{\tilde\Lambda_{x_0}^{1/3}u}_{L^2}^2}.
\end{eqnarray*}
Then we make use of  Cauchy-Schwarz inequality and the interpolation
inequality that
\[
 \forall~\tilde \eps>0,\quad\norm{\tilde\Lambda_{x_0}^{1/3}u}_{L^2}^2
 \leq
 \tilde\eps\norm{\tilde\Lambda_{x_0}^{2/3}u}_{L^2}^2+
 C_{\tilde\eps}\norm{u}_{L^2}^2,
\]
to obtain the desired estimate \reff{L2}.

{\bf \emph{Step C}.}~Now the equality
\[
  {\rm Re} \comi{P_{x_0}u,~\ell_{x_0}u}_{L^2}={\rm Re}
  \comi{iQ_{x_0}u,~\ell_{x_0}u}_{L^2}+
  {\rm Re} \sum_{j=1}^n\comi{L_j^*L_ju,~\ell_{x_0}u}_{L^2}
\]
gives
\begin{eqnarray*}
  {\rm Re}
  \comi{iQ_{x_0}u,~\ell_{x_0}u}_{L^2}\leq
  \norm{P_{x_0}u}_{L^2}^2+\sum_{j=1}^n\norm{L_j^*L_ju}_{L^2}^2
  +\norm{\ell_{x_0}u}_{L^2}^2\\
  \leq\eps\norm{\tilde\Lambda_{x_0}^{2/3}u}_{L^2}^2+
  C_\eps
  \set{\norm{P_{x_0}u}_{L^2}^2+\norm{u}_{L^2}^2},
\end{eqnarray*}
the last inequality following from \reff{LJ} and \reff{L2}. This
along with \reff{Qell} gives at once
\begin{eqnarray*}
  \norm{\tilde\Lambda_{x_0}^{2/3}u}_{L^2}^2&\leq&
  {\rm Re}\comi{iQ_{x_0}u,~\ell_{x_0}u}_{L^2}+\norm{u}_{L^2}^2\\
  &\leq&\eps\norm{\tilde\Lambda_{x_0}^{2/3}u}_{L^2}^2+
  C_\eps
  \set{\norm{P_{x_0}u}_{L^2}^2+\norm{u}_{L^2}^2}.
\end{eqnarray*}
Letting $\eps$ small enough such that $\eps\leq {1\over2},$ we
obtain the desired upper bound of the  term on the left of
\reff{ES2010}, completing the proof of Proposition \ref{++PA+}. \qed

\bigskip
\noindent{\bf Acknowledgements.} The author would like to thank
Professor N. Lerner for his fruitful discussions and  help during
the preparation of this paper.

%\begin{acknowledgements}
%The author would like to thank Professor N. Lerner for fruitful
%discussions, and his help during the preparation of this paper.
%\end{acknowledgements}

\end{document}